\documentclass[11pt,openbib]{article}

\usepackage{latexsym}
\usepackage{amsmath}
\usepackage{amssymb}
\usepackage{amsfonts}
\usepackage{mathrsfs}

\usepackage[breaklinks,
  colorlinks=true,
  citecolor=cyan,
  linkcolor=magenta,
  urlcolor=blue,
  hyperfootnotes=false]{hyperref}

\numberwithin{equation}{section}
\newtheorem{thrm}{Theorem}[section]

\newtheorem{defn}{Definition}[section]
\newtheorem{lemm}{Lemma}[section]
\newtheorem{prop}{Proposition}[section]
\newtheorem{rmrk}{Remark}[section]

\newenvironment{rmk}{{\bf Remark:}}{\hfill{$\Box$}}
\newenvironment{prf}{{\em Proof:}}{\hfill{$\Box$}}

\newcommand{\alf}{\alpha}
\newcommand{\apx}{\approx}

\newcommand{\C}{\mathscr{C}}

\newcommand{\done}{\hfill{$\Box$}}

\newcommand{\dt}[1]{\frac{d}{dt}{#1}}

\newcommand{\Dlt}{\Delta}

\newcommand{\field}[1]{\mathbb{#1}}

\newcommand{\Gmm}{\Gamma}
\newcommand{\gs}{\geqslant}

\newcommand{\hf}{\frac{1}{2}}

\newcommand{\lan}{\langle}

\newcommand{\lf}{\left}

\newcommand{\lra}[2]{\left\langle #1,\,#2\right\rangle}
\newcommand{\ls}{\leqslant}
\newcommand{\lsp}{\preccurlyeq}

\newcommand{\nbl}{\nabla}

\newcommand{\ol}{\overline}

\newcommand{\prt}{\partial}
\newcommand{\pdt}[1]{\frac{\partial#1}{\partial t}}
\newcommand{\pdx}[1]{\frac{\partial#1}{\partial x}}
\newcommand{\pdz}[1]{\frac{\partial#1}{\partial z}}

\newcommand{\Px}{\prt_x}
\newcommand{\Pz}{\prt_z}

\newcommand{\R}{\field{R}}
\newcommand{\ra}{\rightarrow}
\newcommand{\ran}{\rangle}
\newcommand{\rau}{\rightharpoonup}

\newcommand{\rt}{\right}

\newcommand{\sm}{\setminus}

\newcommand{\tld}{\widetilde}
\newcommand{\tu}{\tilde{u}}
\newcommand{\tv}{\tilde{v}}
\newcommand{\tq}{\tilde{q}}
\newcommand{\tw}{\tilde{w}}
\newcommand{\tth}{\tilde{\tht}}

\newcommand{\tht}{\theta}

\newcommand{\veps}{\varepsilon}
\newcommand{\vfi}{\varphi}

\newcommand{\x}{\times}
\newcommand{\zt}{\zeta}

\begin{document}
\pagestyle{myheadings}

\title{\bf Uniqueness of Some Weak Solutions\\
 for 2D Viscous Primitive Equations}

\author{ Ning Ju \footnote{
	Department of Mathematics, Oklahoma State University,
	401 Mathematical Sciences, Stillwater, OK 74078, USA.
	Email: {\tt ning.ju@okstate.edu},
	}}
\date{}

\maketitle

\begin{abstract}

 First, a new sufficient condition for uniqueness of weak solutions is proved
 for the system of 2D viscous Primitive Equations. Second, global existence
 and uniqueness are established for several classes of weak solutions with
 partial initial regularity, including but not limited to those weak solutions
 with initial horizontal regularity, rather than vertical regularity. Our
 results and analyses for the problem with {\em physical boundary conditions}
 can be extended to those with other typical boundary conditions. Most of the
 results were not available before, even for the periodic case.\\

\noindent
{\bf Keywords:} viscous Primitive Equations, existence, uniqueness.\\
{\bf MSC:} 35A01, 35A02, 35B40, 35Q10, 35Q35, 35Q86.

\end{abstract}
\tableofcontents
\maketitle

\markboth{Viscous Primitive Equations}
{\hspace{.2in} N. Ju \hspace{.5in} Viscous Primitive Equations}

\indent
\baselineskip 0.58cm

\section{Introduction}
\label{s:in}

 Consider as in \cite{p.t.z:09} the system of 2D viscous Primitive Equations
 (PE) for {\em three} dimensional Geophysical Fluid Dynamics in the {\em two}
 dimensional spatial domain:
\[ D :=\lf\{(x,z)\in\R^2\ \big|\ 0<x<1, -h<z<0 \rt\}, \]
 where $h$ is a positive constant.

\noindent
 {\em Horizontal momentum equations:}
\begin{align*}
 \pdt{u} + (u, w)\cdot \nbl u =& -\pdx{p}  +v   + \Dlt u,\\
 \pdt{v} + (u,w) \cdot \nbl v =& -u + \Dlt v.
\end{align*}
 {\em Hydrostatic balance:}
\begin{equation*}
 \pdz{p} + \tht =0.
\end{equation*}
 {\em Continuity equation:}
\begin{equation*}
  \pdx{u} + \pdz{w} = 0.
\end{equation*}
 {\em Heat equation:}
\begin{equation*}
  \pdt{\tht} + (u, w)\cdot \nbl\tht = \Dlt \tht +Q.
\end{equation*}
 In the above equations, the gradient $\nbl$ and the Laplacian $\Dlt$ are
 defined as: 
\[ \nbl := (\Px, \Pz)=(\prt_1,\prt_2), \quad
   \Dlt:=\Px^2+\Pz^2=\prt_1^2+\prt_2^2. \]
 The unknowns in the above system of 2D viscous PE are the fluid velocity
 $(u,v, w)\in\R^3 $,  the pressure $p$ and the temperature $\tht$.
 The heat source $Q$ is given. For issues concerned in this article and for
 simplicity of presentation, $Q$ is assumed to be independent of $t$.
 Upon minor modifications, all the results obtained in this article can be
 extended to the case for time-dependent $Q$ under suitable assumptions for $Q$.
 Some of the coefficients in the above system are already simplified for
 conciseness of presentation. In particular, viscosity constant, diffusivity
 constant and Coriolis rotational frequency from $\beta$-plane approximation
 are set as $1$. The effect of salinity is omitted for simplicity of
 presentation. All these simplifications lose no mathematical generality.

 The boundary $\prt D$ of $D$ is partitioned into three parts
 $\Gmm_i \cup \Gmm_b \cup \Gmm_l$, where
\[ \Gmm_i := \{(x, 0) \in \ol{D} \}, \quad
 \Gmm_b := \{(x, -h) \in \ol{D} \}, \quad
 \Gmm_l := \{(x,z)\in\ol{D} \big| x=0,1 \}. \]
 The following boundary conditions of the PEs are used:
\begin{equation*}
\begin{aligned}
  \text{on}\ \Gmm_i \mbox{\,:}& \quad 
  u_z + \alf_1 u = v_z + \alf_2 v= w= \tht_z + \alf_3\tht=0,\\
  \text{on}\ \Gmm_b \mbox{\,:}& \quad u=v=w=\tht_z=0,\\
  \text{on}\ \Gmm_l \mbox{\,:}& \quad u=v=\tht_x=0,
\end{aligned}
\end{equation*}
 where $\alf_i\gs0$ for $i=0,1,2$ and $u_z$, $v_z$, $\tht_x$ and $\tht_z$ are
 the corresponding partial derivatives of $u$, $v$ and $\tht$. For convenience
 of reference, this set of boundary conditions is called {\em physical boundary
 conditions} here, following \cite{e.g:13}. See \cite{l.t.w:92} and
 \cite{p.t.z:09} for geophysical background of these boundary conditions.

 Integrating the continuity equation and the hydrostatic balance equation and
 using the boundary condition $w(x,0,t)=0$, one can express $w$ and $p$ as:
\begin{align}
\label{e:w}
 w(x, z, t) &=  \int_z^0 u_x(x, \zt, t)d\zt,\\
\label{e:p}
 p(x, z, t) &=  q(x, t) +\int^0_z \tht(x, \zt, t)d\zt,
\end{align}
 where $q(x,t)=p(x,0,t)$ is the pressure on $\Gmm_i$.
 Eliminating $w$ and $p$ from the previous system of 2D viscous PE results
 in the following equivalent formulation:
\begin{align}
 u_t -\Dlt u + uu_x &
 +\left(\int_z^0 u_x(x,\zt,t)d\zt\right)u_z \nonumber \\
\label{e:u}
 & + q_x + \int_z^0 \tht_x(x,\zt,t)d\zt = v, \\
\label{e:v}
 v_t -\Dlt v + uv_x &
 +\left(\int_z^0 u_x(x,\zt,t)d\zt\right)v_z = -u, \\
\label{e:t}
 \tht_t -\Dlt\tht + u\tht_x 
 & +\left(\int_z^0 u_x(x,\zt,t)d\zt\right) \tht_z
 =Q;
\end{align}
 with the following boundary conditions:
\begin{align}
\label{e:bc.u}
  u_z+\alf_1 u \big|_{z=0}= u\big|_{z=-h} =u\big|_{x=0,1} =0,\\
\label{e:bc.v}
  v_z +\alf_2 v \big|_{z=0}= v\big|_{z=-h} =v\big|_{x=0,1} =0,\\
\label{e:bc.w}
 \int_{-h}^0 u_x(x,\zt,t)d\zt =0,\\
\label{e:bc.t}
 \tht_z+\alf_3\tht \big|_{z=0} = \tht_z\big|_{z=-h} =\tht_x\big|_{x=0,1} =0.
\end{align}
 The above system of 2D viscous PE will be solved with suitable initial
 conditions:
\begin{equation}
\label{e:ic}
 u(x, z, 0) = u_0(x, z), \quad v(x, z, 0) = v_0(x, z), \quad
 \tht(x, z, 0) = \tht_0(x, z).
\end{equation}
 Notice that the boundary condition $w(x,0,t)=0$ is already embedded in the
 expression \eqref{e:w}. The other boundary condition for $w$ is given in
 \eqref{e:bc.w}. It follows from \eqref{e:bc.u} and \eqref{e:bc.w} that
\begin{align}
\label{e:bc.u.2}
\int_{-h}^0 u(x,z,t)\ dz=0.
\end{align}

 The mathematical framework of the viscous primitive equations for large scale
 ocean flow in {\em three} dimensional spatial domain (3D viscous PE) was
 formulated in \cite{l.t.w:92}, where the notions of weak and strong solutions
 were defined and existence of weak solutions was proved.
 Existence of strong solutions {\em local in time} and their uniqueness were
 proved in \cite{g.m.r:01} and \cite{t.z:04}. Existence of strong solutions
 {\em global in time} was proved in \cite{c.t:07} and \cite{kob:07} for the
 case when $u$ and $v$ satisfy {\em Neumann} boundary condition at the bottom.
 Existence of strong solutions {\em global in time} was proved in \cite{k.z:07}
 for the case when $u$ and $v$ satisfy {\em physical} boundary conditions. See
 also the results in \cite{h.k:16}.
 Uniform boundedness in $H^1$ of strong solutions global in time was proved in
 \cite{j:07} and \cite{k.z:08}. Uniform boundedness in $H^2$ of $H^2$ solutions
 global in time was proved in \cite{j:17} and \cite{j.t:15}. Global uniform
 boundedness in $H^m$ ($m\gs2$) of $H^m$ solutions was recently proved in
 \cite{j:17.10a}.
 
 One of the outstanding unresolved mathematical problems for 3D viscous PE is
 about uniqueness of weak solutions. Global existence and uniqueness of $z$-weak
 solutions to 3D viscous PE were proved in \cite{t-m:10} for initial data in
 $L^6$. A ``$z$-weak'' solution is a weak solution $(u, v,\tht)$ such that
\[ (u_z, v_z,\tht_z)\in L^\infty(0,T; (L^2)^3)\cap L^2(0,T; (H^1)^3). \]
 Recently, global existence and uniqueness of $z$-weak solutions was proved in
 \cite{j:17}. Uniqueness of weak solutions was proved in \cite{k.p.r.z:14} for
 {\em continuous} initial data as well. See also \cite{l.t:17} for a result on
 uniqueness of weak solutions for a class of discontinuous initial data.

 Indeed, the problem of uniqueness of weak solution is still open even for 2D
 viscous PE. In \cite{b.g.m.r:03b}, existence and uniqueness of $z$-weak
 solutions, in the name of ``{\em weak vorticity solutions}'', for 2D
 hydrostatic Navier-Stokes equations (2D hNSE) were proved for the case when $u$
 satisfies a {\em Robin} type {\em friction} boundary condition at the bottom
 of spatial domain. The system of 2D hNSE is somewhat simplified from that of
 2D PE \eqref{e:u}-\eqref{e:t}. It includes only $u$, $w$ and $q$ as the unknown
 variables without $v$ or $\tht$.
 Existence and uniqueness of $z$-weak solution $u$ for 2D hNSE were also proved
 in \cite{b.k.l:04} for the case with {\em Dirichlet} boundary condition at the
 bottom. For the case when the spatial domain is a rectangle, existence and
 uniqueness of the weak solution $u$ of the 2D hNSE were proved in
 \cite{b.k.l:04} with even less demanding regularity:
\[ u \in L^\infty(0, T; L^2_xH^\hf_z)\cap L^2(0,T; H^1_xH^{\frac{3}{2}}_z), \]
 for both the cases of {\em Dirichlet} and {\em Neumann} boundary conditions at
 the bottom. See \cite{b.k.l:04} for notational details. Finally, we mention
 that existence and uniqueness of $z$-weak solution $(u,v,\tht)$ for 2D viscous
 PE were proved in \cite{p:07} for the case with {\em periodic} boundary
 conditions on $(u,v,\tht)$.

 This paper will focus on the problem of uniqueness of weak solutions and
 uniform boundedness of norms of partial regularity for given initial partial
 regularity of weak solutions to the 2D viscous PE. It studies existence and
 uniqueness of weak solutions of the system \eqref{e:u}-\eqref{e:t} in $D$
 under {\em physical} boundary conditions \eqref{e:bc.u}-\eqref{e:bc.t}.

 First, a new sufficient condition for uniqueness of weak solutions is proved.
 Second, global existence of several classes of weak solutions with initial
 partial regularity is also proved. Finally, as an application of our new
 sufficient condition, uniqueness of these classes of weak solutions with
 initial {\em partial} regularity is also proved. These results are valid
 as well for other typical boundary conditions for 2D viscous PEs, since our
 proofs can be easily extended to those cases.

 To present our analysis in complete details, we first give the definition of
 a weak solution carefully and then prove several important results about
 properties of weak solutions and strong solutions of 2D viscous PEs. These
 results are included in Theorems~\ref{t:ws}-\ref{t:ee.ws}. Closely related
 important discussions are also presented in Remarks~\ref{r:ws.def}-\ref{r:lc}.
 This section, Section~\ref{s:ws}, shares some similarity with \cite{l.t:17}
 in terms of strategy. However, the definition of a weak solution of the
 viscous PE used in this paper is somewhat {\em different} from those used in
 \cite{b.g.m.r:03b}, \cite{l.t:17} and \cite{l.t.w:92} for viscous PE. The
 set of boundary conditions used in this paper is also {\em different} from
 those used in \cite{b.g.m.r:03b} and \cite{l.t:17}. Hence, many of our
 detailed arguments and ideas of the proofs are also {\em different} from those
 of \cite{l.t:17}. Therefore, we choose to present the full proofs of all these
 results, which will provide fundamental technical support for our analysis in
 the rest sections of this paper. Some of these results might also be new in
 the presented forms.

 The main results for existence and uniqueness of weak solutions with initial
 partial regularity to be presented are Theorem~\ref{t:sc}, Theorem~\ref{t:uz},
 Theorem~\ref{t:ux} and Theorem~\ref{t:uniq}. Except that the results for
 $z$-weak solutions were proved in \cite{p:07} for periodic case, all of our
 main results were not previously known even for periodic case. Especially,
 global existence and uniqueness are proved for weak solutions with initial
 partial regularity in horizontal direction, rather than vertical direction.
 The same is proved for several other mixed cases as well. The main ideas of
 our analysis are careful manipulations of anisotropic inequalities in Sobolev
 spaces.

 The rest of this article is organized as follows:

 In Section~\ref{s:pr}, we give the notations and some definitions, briefly
 review the background and recall some facts and known results which are
 important for later analysis.
 In Section~\ref{s:ws}, we first give the definition of weak solutions and
 strong solutions of 2D viscous PE with physical boundary conditions. Then,
 prove Theorem~\ref{t:ws}-\ref{t:ee.ws} about important properties of weak
 solutions and strong solutions of 2D viscous PE. These theorems will provide
 helpful technical support for our analysis in the sequel.
 In Section~\ref{s:sc}, we prove Theorem~\ref{t:sc}. It gives a new sufficient
 condition for uniqueness of weak solutions and generalizes an already known
 one. It is a crucial observation which will be used later in proving our new
 uniqueness results.
 In Section~\ref{s:ex}, we prove Theorem~\ref{t:uz} and Theorem~\ref{t:ux}.
 These global existence results, for weak solutions in specific partial
 regularity classes, not only give corresponding  global-in-time uniform bounds
 and absorbing sets, but also prepare our proof of uniqueness in next section.
 In Section~\ref{s:un}, we prove Theorem~\ref{t:uniq} for uniqueness of the
 weak solutions with the specific initial partial regularity.

\section{Preliminaries}
\label{s:pr}

\noindent
 In this paper, $C$ denotes a positive absolute constant, the value of which
 might vary from line to line. Similarly, $C_\veps$ denotes a positive constant
 depending on $\veps>0$, the value of which may also vary at different
 occurrence. The following notations are used for real numbers $A$ and $B$:
\[  A \lsp B \quad \text{ if and only if }\quad  A\ls C\cdot B,\]
 and
\[  A \apx B \quad \text{ if and only if }\quad  c\cdot A\ls B \ls C\cdot A, \]
 for some positive constants $c$ and $C$ independent of $A$ and $B$.

 Denote by $L^r(D)$, $L^r((0,1))$ and $L^r((-h,0))$ ($1\ls r <+\infty$) the
 classic Lebesgue $L^r$ spaces with the norm:
 \begin{equation*}
 \|\vfi\|_r = \left\{ \begin{array}{ll}
   \left( \int_D|\vfi(x,z)|^r\ dxdz\right)^\frac{1}{r},
	& \forall \vfi\in L^r(D); \\
  \left( \int_0^1|\vfi(x)|^r\ dx \right)^\frac{1}{r},
	& \forall \vfi\in L^r((0,1));\\
  \left( \int_{-h}^0|\vfi(z)|^r\ dz \right)^\frac{1}{r},
	& \forall \vfi\in L^r((-h,0)).
  \end{array} \right.
 \end{equation*}
 Standard modification is used when $r=\infty$. When there is no confusion,
 index $r=2$ is omitted:
\[ \|\vfi\| := \|\vfi\|_2. \]
 Denote by $H^m(D)$ ($m\gs 1$) the classic Sobolev spaces for square-integrable
 functions on $D$ with square-integrable weak derivatives up to order $m$.
 Domains of the functions spaces will be omitted from notations without
 confusion.

 Some {\em anisotropic} Lebesgue spaces and Sobolev spaces will be used.
 For example, for $r,s\in[1,\infty)$, $L^r_x(L^s_z)$ denotes the standard
 function space of (classes of) Lebesgue measurable functions on $D$ such that
\begin{equation*}
\begin{split}
\|\vfi\|_{L^r_x(L^s_z)}
 &:=\lf(\int_0^1 \|\vfi(x,\cdot)\|_{L^s_z}^r dx\rt)^{\frac{1}{r}}\\
 &:=\lf[\int_0^1
  \lf(\int_{-h}^0 |\vfi(x,z)|^s dz\rt)^{\frac{r}{s}} dx \rt]^{\frac{1}{r}}
 <\infty,
\end{split}
\end{equation*}
 with stand modifications when $r$ or $s$ is $\infty$.

 We will also use $C_B([\alf,\infty))$ to denote the set of uniformly bounded
 functions in $C([\alf,\infty))$, for an interval $[\alf,\infty)\subset\R$.

 Following \cite{p.t.z:09}, the following function spaces are defined:
\begin{equation*}
\label{e:HnV}
 H := H_1\x H_2\x H_3,\quad V:= V_1\x V_2\x V_3,
\end{equation*}
 where
\begin{equation*}
\begin{split}
 H_1 &:= \lf\{ \vfi \in L^2(D)\ \Big|\ \int_{-h}^0 \vfi(x,z)dz =0,
 \vfi\big|_{\Gmm_l}=0 \rt\}, \\
 V_1 &:= \lf\{ \vfi\in H^1(D)\ \Big|\ \int_{-h}^0 \vfi(x,z)dz =0,\
  \vfi\big|_{\Gmm_l\cup\Gmm_b}=0 \rt \} , \\
 H_2 &:= L^2(D),\
 V_2 :=\lf\{\vfi\in H^1(D)\ \Big|\ \vfi\big|_{\Gmm_l\cup\Gmm_b}=0 \rt\},
\end{split}
\end{equation*}
\[ H_3 := L^2(D),\quad V_3 := H^1(D), \quad \text{ if } \alf_3>0,\]
 and 
\[  H_3 := \lf\{ \vfi \in L^2(D)\ \Big|\ \int_D \vfi = 0 \rt\}, \
    V_3 := H_3 \cap H^1(D), \quad \text{ if } \alf_3=0.\]
 Define the bilinear form: $a_i: V_i\times V_i \ra \R$ for $i=1,2,3$, such that
\begin{equation*}
\label{e:ai}
 a_i(\phi, \vfi) = \int_D \nbl\phi \cdot \nbl\vfi\ dxdz
 + \alf_i \int_0^1 \phi(x,0)\vfi(x,0)\ dx,
\end{equation*}
 and the corresponding linear operator $A_i : V_i \mapsto V_i'$, such that
\begin{equation*}
\label{e:Ai}
 \lra{A_iv}{\vfi}=a_i(v, \vfi), \quad \forall v, \vfi \in V_i,
\end{equation*}
 where $V_i'$ is the dual space of $V_i$ and $\lra{\cdot}{\cdot}$ denotes
 scalar product between $V_i'$ and $V_i$ {\em and} the inner product in $H_i$.
 Define:
 \[ D(A_i) =\lf\{ \phi \in V_i\ \big|\ A_i\phi\in H_i \rt\},\quad i=1,2,3. \]
 Define: $A: V\mapsto V'$ as $A(u,v,\tht) := (A_1u,A_2v,A_3\tht)$, for
 $(u,v,\tht)\in V$. Then,
\[ D(A) = D(A_1)\x D(A_2)\x D(A_3). \]
 Since $A_i^{-1}$ is a self-adjoint compact operator in $H_i$, by the classic
 spectral theory, the power $A_i^s$ can be defined for any $s\in\R$. Then,
\[ D(A_i)' = D(A_i^{-1})\]
 is the dual space of $D(A_i)$ and
\[ V_i = D(A_i^\hf), \quad V_i' = D(A_i^{-\hf}).\] 
 Moreover,
\[ D(A_i) \subset V_i \subset H_i \subset V_i' \subset D(A_i)', \]
 and the embeddings above are all continuous and compact and each space
 above is dense in the one following it. Define the norm of $V_i$ as
 \[ \|\vfi\|_{V_i}^2 = a_i(\vfi,\vfi) = \lra{A_i\vfi}{\vfi}
 = \lan A_i^\hf\vfi, A_i^\hf\vfi \ran, \quad i=1,2,3. \]
 Then, for $\vfi\in V_i$ and $i=1,2,3$,
\begin{equation*}
\label{e:ponc.nmeqv}
 \|\vfi\| \lsp \|\vfi\|_{V_i} \apx \|\vfi\|_{H^1}.
\end{equation*}
 By the above discussion and elliptic regularity for linear 3D stationary PE
 (see e.g. \cite{t.z:04}), we also have for $\vfi\in D(A_i)$ and $i=1,2,3$,
\begin{equation*}
\label{e:h2.nmeqv}
 \|\vfi\|_{V_i} \lsp \|A_i\vfi\| \apx \|\vfi\|_{H^2}.
\end{equation*}
 Next, we introduce the following anisotropic estimate which will be very useful
 for our later discussion:
\begin{lemm}
\label{l:tri}
 Let $\psi,\psi_x,\phi,\phi_z,\vfi \in L^2(D)$. Then,
\[ \lf|\int_D \psi\phi\vfi \rt|
 \lsp \|\psi\|^\hf(\|\psi\|+\|\psi_x\|)^\hf
     \|\phi\|^\hf(\|\phi\|+\|\phi_z\|)^\hf
     \|\vfi\|. \]
\end{lemm}

\begin{prf}
\begin{equation*}
\begin{split}
 \lf|\int_D \psi\phi\vfi\rt|
\ls&  \int_0^1 \|\phi\|_{L^\infty_z}\|\psi\|_{L^2_z}\|\vfi\|_{L^2_z}\ dx\\
\lsp& \int_0^1 \|\phi\|_{L^2_z}^\hf(\|\phi\|_{L^2_z}+\|\phi_z\|_{L^2_z})^\hf
      \|\|\psi\|_{L^2_z}\|\vfi\|_{L^2_z} dx\\
\lsp& \|\psi\|_{L^\infty_x(L^2_z)}\int_0^1 \|\phi\|_{L^2_z}^\hf
      (\|\phi\|_{L^2_z}+\|\phi_z\|_{L^2_z})^\hf\|\vfi\|_{L^2_z}\ dx\\
\lsp& \|\psi\|_{L^2_z(L^\infty_x)}\|\phi\|^\hf (\|\phi\|+\|\phi_z\|)^\hf\|\vfi\|\\
\lsp& \|\psi\|^\hf(\|\psi\|+\|\psi_x\|)^\hf
      \|\phi\|^\hf(\|\phi\|+\|\phi_z\|)^\hf\|\vfi\|.
\end{split}
\end{equation*}
 Notice that, in the last step of the above derivations, we have used the
 following estimate:
\begin{equation*}
\begin{split}
 \|\psi\|_{L^2_z(L^\infty_x)}^2
=&\int_{-h}^0 \|\psi\|_{L^\infty_x}^2 dz\\
\lsp& \int_{-h}^0 \|\psi\|_{L^2_x}(\|\psi\|_{L^2_x}+\|\psi_x\|_{L^2_x})\\
\lsp& \|\psi\| (\|\psi\|+\|\psi_x\|).
\end{split}
\end{equation*}

\end{prf}

\noindent
 An immediate application of Lemma~\ref{l:tri} is the following result:
\begin{lemm}
\label{l:c.nl} The following statements are valid:
\begin{enumerate}
\item[{\em (a)}]
 Suppose $(u,v,\tht)\in V$ and $w$ is given by \eqref{e:w}. Then,
 for $\vfi\in V_1$,
\begin{align}
\label{e:u.ux.f}
 \lf|\lra{u u_x}{\vfi}\rt| \lsp& \|u\|\|u_x\|^\hf\|u_z\|^\hf\|\vfi_x\|, \\
\label{e:w.uz.f}
 \lf|\lra{w u_z}{\vfi}\rt| \lsp& \|u\| \|u_x\|^\hf\|u_z\|^\hf\|\vfi_x\|
 + \|u\|^\hf\|u_x\|^{\frac{3}{2}}\|\vfi_z\|;
\end{align}
 for $\vfi\in V_2$, with $i, j=1,2$, $i'=3-i$ and $j'=3-j$,
\begin{align}
\label{e:u.vx.f}
 \lf|\lra{u v_x}{\vfi}\rt| \lsp
 \|u_x\| &\|v\|^\hf\|\prt_iv\|^\hf\|\vfi\|^\hf\|\prt_{i'}\vfi\|^\hf
 \nonumber\\
 +& \|u\|^\hf\|\prt_j u\|^\hf \|v\|^\hf\|\prt_{j'}v\|^\hf\|\vfi_x\|, \\
 \lf|\lra{w v_z}{\vfi}\rt|
\lsp \|u_x\| &\|v\|^\hf\|\prt_iv\|^\hf\|\vfi\|^\hf\|\prt_{i'}\vfi\|^\hf
\nonumber\\
\label{e:w.vz.f}
 &+\|u_x\|\|v\|^\hf\|v_x\|^\hf\|\vfi_z\|;
\end{align}
 for $\vfi\in V_3$, with $i, j=1,2$, $i'=3-i$ and $j'=3-j$,
\begin{align}
\label{e:u.tx.f}
 \lf|\lra{u \tht_x}{\vfi}\rt| \lsp
 \|u_x\| &\|\tht\|^\hf\|\prt_i\tht\|^\hf\|\vfi\|^\hf\|\prt_{i'}\vfi\|^\hf
 \nonumber\\
 +& \|u\|^\hf\|\prt_j u\|^\hf \|\tht\|^\hf\|\prt_{j'}\tht\|^\hf\|\vfi_x\|, \\
 \lf|\lra{w \tht_z}{\vfi}\rt| \lsp
 \|u_x\| &\|\tht\|^\hf\|\prt_i\tht\|^\hf\|\vfi\|^\hf\|\prt_{i'}\vfi\|^\hf
\nonumber\\
\label{e:w.tz.f}
 &+\|u_x\|\|\tht\|^\hf\|\tht_x\|^\hf\|\vfi_z\|.
\end{align}
 Therefore,
\begin{equation}
\label{e:V*}
 uu_x, wu_z\in V_1',\quad uv_x, wv_z\in V_2',\quad u\tht_x, w\tht_z\in V_3'.
\end{equation}

\item[{\em (b)}]
 Suppose $(u,v,\tht)\in L^\infty(0,T; H)\cap L^2(0,T; V)$. Then,
\begin{equation*}
\begin{split}
 uu_x\in L^2(0, T; V_1'), \quad & wu_z \in L^{\frac{4}{3}}(0, T; V_1'),\\
 uv_x, wv_z \in L^{\frac{4}{3}}(0, T; V_2'), \quad & 
 u\tht_x, w\tht_z \in L^{\frac{4}{3}}(0, T; V_3').
\end{split}
\end{equation*}

\end{enumerate}
\end{lemm}

\begin{prf}
 For $u,\vfi\in V_1$, Lemma~\ref{l:tri} yields
\[ \lf|\lra{u u_x}{\vfi}\rt|
 \lsp \|u\|^\hf \|u_x\|^{\frac{3}{2}}\|\vfi\|^\hf\|\vfi_z\|^\hf, \quad
 \lf|\lra{u^2}{\vfi_x}\rt| \lsp \|u\|\|u_x\|^\hf\|u_z\|^\hf\|\vfi_x\|. \]
 Then, a density argument using Lemma~\ref{l:tri} again, along with the above
 two inequalities, proves
\[ \lra{uu_x}{\vfi} = -\hf \lra{u^2}{\vfi_x},\]
 from which \eqref{e:u.ux.f} follows. Similarly, Lemma~\ref{l:tri} yields
\[ \lf|\lra{w u_z}{\vfi}\rt|
  \lsp \|w\|^\hf\|w_z\|^\hf\|\vfi\|^\hf\|\vfi_x\|^\hf\|u_z\|\\
 \ls \|u_x\| \|u_z\| \|\vfi\|^\hf\|\vfi_x\|^\hf, \]
 and
\begin{equation}
\label{e:wuf}
 \lf|\lra{w u}{\vfi_z}\rt|
 \lsp \|u\|^\hf\|u_x\|^{\frac{3}{2}}\|\vfi_z\|.
\end{equation}
 The above two inequalities plus \eqref{e:u.ux.f}, along with Lemma~\ref{l:tri},
 then imply via a density argument that
\[ \lra{wu_z}{\vfi} = \lra{u_xu}{\vfi}- \lra{w u}{\vfi_z}, \]
 from which, we immediately prove \eqref{e:w.uz.f} by \eqref{e:u.ux.f} and
 \eqref{e:wuf}.

 Similarly, we can prove \eqref{e:u.vx.f}-\eqref{e:w.tz.f}. Then, it is easy
 to prove \eqref{e:V*} and part (b) using \eqref{e:u.ux.f}-\eqref{e:w.tz.f}.

\end{prf}

\section{Weak Solutions and Strong Solutions}
\label{s:ws}

 In this section, some important properties about weak solutions and strong
 solutions will be discussed. These will provide important technical support
 in the proofs of the main results of this paper to be presented in the next
 few sections.

 The following definitions of weak and strong solutions of the initial boundary
 value problem \eqref{e:u}-\eqref{e:ic} for the 2D viscous PEs will be used in
 this paper:
\begin{defn}
\label{d:soln}
{\em Suppose $Q \in L^2(D)$, $(u_0, v_0, \tht_0)\in H$ and $T>0$. The triple
 $(u, v, \tht)$ is called a {\em weak solution} of the viscous PEs
 \eqref{e:u}-\eqref{e:ic} on the time interval $(0, T)$ if it satisfies
 \eqref{e:u}-\eqref{e:t} in weak sense, 
 that is, if
\begin{equation}
\label{e:reg.w}
 (u,v,\tht)\in L^\infty(0, T;H)\cap L^2(0,T; V),
\end{equation}
 satisfies the follow equations in the sense of distribution on $(0,T)$:
\begin{align}
\label{e:u.w} 
\lra{u_t}{\vfi} + a_1(u,\vfi)
 + \lra{(u,w)\cdot\nbl u -v+\int_z^0\tht_x}{\vfi}
 &= 0, \quad \forall \vfi \in V_1,\\
\label{e:v.w} 
\lra{v_t}{\vfi} + a_2(v,\vfi) + \lra{(u,w)\cdot\nbl v -u}{\vfi}
 &=0, \quad \forall \vfi \in V_2,\\
\label{e:t.w}
 \lra{\tht_t}{\vfi} + a_3(\tht,\vfi)
 + \lra{(u,w)\cdot\nbl\tht- Q}{\vfi} &=0,
 \quad \forall \vfi \in V_3,
\end{align}
 where $w$ is given by \eqref{e:w} in weak sense. Moreover,
\begin{equation}
\label{e:ic.w}
 \lim_{t\ra 0^+}(u(t),v(t),\tht(t)) = (u_0,v_0,\tht_0),
\end{equation}
 in weak topology of $H$, and the following {\em energy inequalities} are
 satisfied for {\em almost every} $t_0\in[0, T)$ and {\em almost every}
 $t\in(t_0,T)$:
\begin{equation}
\label{e:ei.u}
\|u(t)\|^2
 + 2\int_{t_0}^{t}\lf(\|u(s)\|_{V_1}^2+\lra{-v+\int_z^0\tht_x}{u}\rt) ds
 \ls \|u(t_0)\|^2,
\end{equation}
\begin{equation}
\label{e:ei.v}
\|v(t)\|^2 +2\int_{t_0}^t\lf(\|v(s)\|_{V_2}^2+\lra{u}{v}\rt) ds
 \ls \|v(t_0)\|^2,
\end{equation}
\begin{equation}
\label{e:ei.t}
 \|\tht(t)\|^2
+ 2\int_{t_0}^t \lf(\|\tht(s)\|_{V_3}^2 -\lra{Q}{\tht(s)}\rt) ds
 \ls \|\tht(t_0)\|^2.
\end{equation}
 Further more, the above energy inequalities \eqref{e:ei.u}-\eqref{e:ei.t} are
 also satisfied for $t_0=0$ and for {\em almost every} $t\in(0, T)$.

 If $(u_0, v_0, \tht_0)\in V$, then $(u, v, \tht)$ is called a
 {\em strong solution} of \eqref{e:u}-\eqref{e:ic} on the time interval
 $[0, T)$ if it satisfies \eqref{e:u.w}-\eqref{e:ic.w} and
\begin{equation}
\label{e:reg.s}
 (u, v, \tht) \in L^\infty(0, T; V) \cap L^2(0, T; D(A)).
\end{equation}
 If $T>0$ in the above can be arbitrarily large, then the corresponding weak
 or strong solution is {\em global}.}  \done 
 
\end{defn}

\begin{rmrk}
\label{r:ws.def}
{\em There are somewhat different ways to define weak solutions of the PE.
 For examples, see \cite{b.g.m.r:03a},\cite{b.g.m.r:03b}, \cite{l.t:17},
 \cite{l.t.w:92}, \cite{p.t.z:09} and \cite{t.z:04}). Especially, to define
 a weak solution of the 3D PE with physical boundary conditions, the domain
 of $\vfi$ in \eqref{e:u.w}-\eqref{e:t.w} was chosen as $D(A_i)$ in
 \cite{l.t.w:92}, \cite{p.t.z:09} and \cite{t.z:04}, for $i=1,2,3$ respectively
 instead of $V_i$. Definition~\ref{d:soln} is formally more restrictive than
 the one given in \cite{l.t.w:92}, \cite{p.t.z:09} and \cite{t.z:04}. However,
 $D(A_i)$ is dense in $V_i$ and, by Lemma~\ref{l:c.nl}, the nonlinear terms of
 the 2D PE are in $V_i'$ for $(u,v,\tht)\in V$ and in $L^{\frac{4}{3}}(0,T;V_i')$
 for $(u,v,\tht)\in L^\infty(0,T;H)\cap L^2(0,T;V)$.}
 Thus, for 2D case, a weak solution defined in \cite{l.t.w:92}, \cite{p.t.z:09}
 and \cite{t.z:04}, if satisfying \eqref{e:ei.u}-\eqref{e:ei.t}, is also a weak
 solution in the sense of Definition~\ref{d:soln}. \done

\end{rmrk}

 We first state and prove the following theorem on some basic properties
 satisfied by every weak solution.

\begin{thrm}
\label{t:ws}
 There exists at least one global weak solution of \eqref{e:u}-\eqref{e:ic} in
 the sense of Definition~\ref{d:soln}. If $(u,v,\tht)$ is a weak solution on
 $(0, T)$, then\footnote{If $T=\infty$, the space $C([0,T])$ in \eqref{e:c.w}
 is replaced by $C_B([0,\infty))$.}
\begin{equation}
\label{e:c.w}
 \lra{u}{\vfi_1}, \lra{v}{\vfi_2},\lra{\tht}{\vfi_3}
 \in C([0,T]), \quad \forall\ (\vfi_1,\vfi_2, \vfi_3)\in H.
\end{equation}
 Moreover, there exists a zero measure set $E\subset (0, \infty)$, such that
\begin{equation}
\label{e:lsc.H}
 \lim_{E^c\ni t\ra 0^+}\|(u(t),v(t),\tht(t))-(u_0,v_0,\tht_0)\|_H =0,
\end{equation}
 where $E^c:=(0,\infty)\sm E$.
\end{thrm}

\begin{prf}
 It is proved in \cite{l.t.w:92} that there exists at least one global weak
 solution (in their sense) for 3D PE and it satisfies energy inequalities
 \eqref{e:ei.u}-\eqref{e:ei.t}. Moreover, any weak solution as defined in
 \cite{l.t.w:92} is weakly continuous from $[0,T]$ into $H$ if $T$ is finite
 and weakly continuous from $[0,T)$ into $H$ if $T=\infty$.
 By Remark~\ref{r:ws.def}, these results imply existence of at least one global
 weak solution $(u,v,\tht)$ of \eqref{e:u}-\eqref{e:ic} in the sense of
 Definition~\ref{d:soln} and that \eqref{e:c.w} is satisfied by {\em any} weak
 solution $(u,v,\tht)$ in the sense of Definition~\ref{d:soln}.

 By Lemma~\ref{l:tri} and Lemma~\ref{l:c.nl}, we can also prove existence of
 at least one global weak solutions $(u,v,\tht)$ of the 2D problem
 \eqref{e:u}-\eqref{e:ic} using Definition~\ref{d:soln} directly, by following
 the standard approach of \cite{t.nse} for Navier-Stokes equations. Moreover,
 we can prove that any weak solution $(u,v,\tht)$ satisfies \eqref{e:c.w}.

 Finally, due to the fact that a weak solution satisfies the energy inequalities
 \eqref{e:ei.u}-\eqref{e:ei.t} for $t_0=0$ and for almost every $t\in [0, T]$
 by Definition~\ref{d:soln}, there exists a zero measure set
 $E\subset(0,\infty)$ such that, for all $t\in E^c$,
\begin{equation}
\label{e:ei.u0}
\|u(t)\|^2
 + 2\int_0^{t}\lf(\|u(s)\|_{V_1}^2+\lra{-v+\int_z^0\tht_x}{u}\rt)ds
 \ls \|u_0\|^2,
\end{equation}
\begin{equation}
\label{e:ei.v0}
\|v(t)\|^2 +2\int_0^t\lf(\|v(s)\|_{V_2}^2+\lra{u}{v}\rt) ds
 \ls \|v_0\|^2,
\end{equation}
\begin{equation}
\label{e:ei.t0}
 \|\tht(t)\|^2
+ 2\int_0^t \lf(\|\tht(s)\|_{V_3}^2 -\lra{Q}{\tht(s)}\rt) ds
 \ls \|\tht_0\|^2.
\end{equation}
 Notice that, by definition, $(u,v,\tht)\in L^2(0, T; V)$. Therefore, taking
 $\limsup$ on both sides of \eqref{e:ei.u0} for $t(\in E^c)\ra0^+$, we have
\[ \limsup_{E^c\ni t\ra0^+}\|u(t)\|^2 \ls \|u_0\|^2.\]
 By \eqref{e:c.w} and weak lower semi-continuity, we also have
\[ \|u_0\|^2 \ls \liminf_{E^c\ni t\ra0^+}\|u(t)\|^2. \]
 Thus,
\[ \lim_{E^c\ni t\ra0^+} \|u(t)\|^2 = \|u_0\|^2.\]
 Hence,
\begin{equation*}
\begin{split}
 \lim_{E^c\ni t\ra0^+} \|u(t)-u_0\|^2
=&\lim_{E^c\ni t\ra0^+}(\|u(t)\|^2 -2\lra{u(t)}{u_0} +\|u_0\|^2)\\
=&\lim_{E^c\ni t\ra0^+}(\|u(t)\|^2 -\|u_0\|^2) =0.
\end{split}
\end{equation*}
 The second equality above is due to weak continuity \eqref{e:c.w}. This weak
 continuity argument was also used in \cite{l.t:17} to prove \eqref{e:lsc.H}.

\end{prf}

 Different proofs of existence of global weak solutions of 3D PE can also be
 found in \cite{b.g.m.r:03a} and \cite{t.z:04}. Different sets of boundary
 conditions and different definitions of weak solution were used in
 \cite{b.g.m.r:03a}, \cite{b.g.m.r:03b} and \cite{l.t:17}.

\begin{rmrk}
\label{r:ss}
 {\em Existence and uniqueness of global strong solution for 3D PE with
 {\em Neumann boundary condition} for $(u,v)$ at bottom was proved in
 \cite{c.t:07}. See also \cite{kob:07} for a different proof of existence of
 global strong solution with the same boundary conditions when initial data is
 in $H^2$. Existence and uniqueness of global strong solution was proved in
 \cite{k.z:07} and \cite{k.z:08} for 3D viscous PE with {\em physical boundary
 condition}. The strong solutions are uniformly bounded in $V$ and a bounded
 absorbing set for the solutions exists in $V$. These results apply to the 2D
 case as well. See also \cite{p.t.z:09} for a direct proof of global existence
 of the strong solution of the 2D viscous PE \eqref{e:u}-\eqref{e:ic}.
 Moreover, following the argument of \cite{j:07} for the case of 3D PE with
 Neumann boundary conditions, we can prove (see \cite{e.g:13}) for the 3D PE
 with physical boundary conditions that, if $(u_0,v_0,\tht_0)\in V$, then
 the strong solution $(u,v,\tht)$ satisfies
\begin{equation}
\label{e:ss.dt}
 (u_t,v_t, \tht_t) \in L^2(0,T;H), \quad\ \forall T>0,
\end{equation}
 and
\begin{equation}
\label{e:ss.c}
  (u(t),v(t),\tht(t)) \in C_B([0, \infty); V).
\end{equation}
} \done
\end{rmrk}
 
 The next theorem shows that energy {\em equalities} are satisfied by every
 strong solution. Therefore, the strong solution is also a weak solution.
\begin{thrm}
\label{t:ee.ss}
 Let $(u,v,\tht)$ be the unique global strong solution of
 \eqref{e:u}-\eqref{e:ic} with $(u_0,v_0,\tht_0)\in V$. Then, for every
 $t_0\in[0,\infty)$ and $t\in(t_0, \infty)$,
\begin{equation}
\label{e:ee.u}
\|u(t)\|^2+2\int_{t_0}^t\lf(\|u(s)\|_{V_1}^2
 + \lra{-v(s)+\int_z^0\tht_x(s)}{u(s)} \rt) ds =\|u(t_0)\|^2,
\end{equation}
\begin{equation}
\label{e:ee.v}
\|v(t)\|^2+2\int_{t_0}^t \lf(\|v(s)\|_{V_2}^2 +\lra{u(s)}{v(s)}\rt)ds
 =\|v(t_0)\|^2,
\end{equation}
\begin{equation}
\label{e:ee.t}
\|\tht(t)\|^2 + \int_{t_0}^t\lf(\|\tht(s)\|_{V_3}^2-\lra{Q}{\tht(s)}\rt) ds =
 \|\tht(t_0)\|^2.
\end{equation}
 Therefore, the strong solution is also a weak solution.
\end{thrm}

\begin{prf}
 By \eqref{e:reg.s}, \eqref{e:ss.dt} and a lemma of Lions and Magenes
 (see \cite{l.m:72} and \cite{t.nse}), we have in the sense of distribution
 on $(0,\infty)$
\begin{equation}
\label{e:lm}
 \frac{d}{dt} \|u\|^2 = 2 \lra{u_t}{u},\
 \frac{d}{dt} \|v\|^2 = 2 \lra{v_t}{v},\
 \frac{d}{dt} \|\tht\|^2 = 2 \lra{\tht_t}{\tht}.
\end{equation}
 Therefore, by \eqref{e:ss.c}, \eqref{e:lm} and Definition~\ref{d:soln},
 we have in {\em classic} sense on $[0, \infty)$,
\begin{equation}
\label{e:ee.du}
 \hf\frac{d}{dt} \|u\|^2 + \|u\|_{V_1}^2 + \lra{\int_z^0\tht_x}{u} = \lra{v}{u},
\end{equation}
\begin{equation}
\label{e:ee.dv}
 \hf\frac{d}{dt} \|v\|^2 +\|v\|_{V_2}^2 + \lra{u}{v} = 0,
\end{equation}
\begin{equation}
\label{e:ee.dt}
 \hf\frac{d}{dt} \|\tht\|^2 + \|\tht\|_{V_3}^2 = \lra{Q}{\tht}.
\end{equation}
 In the above derivation, we have also used the cancellation property for
 $(u,v,\tht)\in V$:
\begin{equation}
\label{e:cancel}
 \lra{(u,w)\cdot\nbl u}{u} = \lra{(u,w)\cdot\nbl v}{v}
 =\lra{(u,w)\cdot\nbl\tht}{\tht}=0,
\end{equation}
 which can be justified by Lemma~\ref{l:tri}, as shown in the proof of
 Lemma~\ref{l:c.nl}.
 Integrating \eqref{e:ee.du}-\eqref{e:ee.dt} finishes the proof.
 
\end{prf}

 {\em Weak-strong uniqueness} was proved in \cite{l.t:17}. That is, a weak
 solution with initial data in the space of strong solutions $V$ is the strong
 solution with the same initial data, and is thus the unique weak (and strong)
 solution. The definition of weak solution used in \cite{l.t:17} is somewhat
 different from Definition~\ref{d:soln}. In the following, we give a completely
 different proof of weak-strong uniqueness result. Our prove is a direct proof
 using Definition~\ref{d:soln} and properties of weak and strong solutions.
 The argument of our proof is also different from that of \cite{t.nse} for a
 related uniqueness result.

\begin{thrm}
\label{t:ws.uniq}
 Let $(u,v,\tht)$ be a weak solution of \eqref{e:u}-\eqref{e:ic} on $(0,T)$ in
 the sense of Definition~\ref{d:soln} with $(u_0,v_0,\tht_0)\in V$. Then,
 $(u,v,\tht)$ is the strong solution of \eqref{e:u}-\eqref{e:ic}.
 Thus, $(u,v,\tht)$ is the unique weak (and strong) solution.
\end{thrm}

\begin{prf}
 Let $(u_1,v_1,\tht_1)$ be the global strong solution and $(u_2,v_2,\tht_2)$ be
 a weak solution on $(0, T)$ for some $T>0$, with the same initial value
 $(u_0,v_0,\tht_0)\in V$.  For convenience of presentation, assume $T$ is
 finite. The case of $T=\infty$ is then an easy consequence.

 Use $w_i$ and $q_i$, for $i=1,2$, to denote corresponding vertical velocity
 and surface pressure. Denote:
\[ (\tu, \tv, \tth, \tw, \tq)
  := (u_1-u_2, v_1-v_2, \tht_1-\tht_2, w_1-w_2, q_1-q_2). \] 
 Notice that, by \eqref{e:ss.dt} and \eqref{e:ss.c}, for the strong solution,
\[ u_1\in C([0,T]; V_1)\cap L^2(0,T; D(A_1)), \quad u_{1,t}\in L^2(0,T; H_1).\]
 For a weak solution $u_2$,
\[ u_2\in L^\infty(0,T; H_1)\cap L^2(0,T; V_1),\quad u_{2,t}\in L^1(0,T;V_1').\]
 Therefore, by standard regularization approximation on $(0,T)$, we can obtain
 sequences of functions $\{u_{1,m}\}_{m=1}^\infty$ and $\{u_{2,m}\}_{m=1}^\infty$
 such that
\[ u_{1,m}\in C^\infty([0,T]; D(A_1)), u_{2,m} \in C^\infty([0,T]; V_1),
 \quad \forall\ m\gs 1, \]
 and as $m\ra\infty$,
\begin{equation}
\label{e:mlfc.prop}
\begin{split}
  u_{1,m} &\ra\ u_1 \quad\text{in}\ L^2_{loc}(0,T; D(A_1)) \text{ and }
 C([0, T]; V_1),\\
 (u_{1,m})_t &\ra\ u_{1,t} \quad\text{in}\ L^2_{loc}(0,T; H_1), \\
 u_{2,m} &\ra\ u_2 \quad \text{in}\ L^2_{loc}(0,T; V_1)\\
 (u_{2,m})_t &\ra\ u_{2,t} \quad\text{in}\ L^1_{loc}(0,T; V_1'),
\end{split}
\end{equation}
 It is obvious that, for any $m\gs 1$,
\begin{equation}
\label{e:prod}
 \frac{d}{dt} \lra{u_{1,m}}{u_{2,m}} = \lra{(u_{1,m})_t}{u_{2,m}}
 + \lra{(u_{2,m})_t}{u_{1,m}}.
\end{equation}
 As $m\ra\infty$, we have by \eqref{e:mlfc.prop} that, in $L^1_{loc}(0,T)$,
\begin{equation}
\label{e:conv}
\begin{split}
 \lra{u_{1,m}}{u_{2,m}} &\ra \lra{u_1}{u_2},\\
 \lra{(u_{1,m})_t}{u_{2,m}} &\ra \lra{u_{1,t}}{u_2},\\
 \lra{(u_{2,m})_t}{u_{1,m}} &\ra \lra{u_{2,t}}{u_1}.
\end{split}
\end{equation}
 These convergences are also valid in the distribution sense. Therefore,
 we can take the limit $m\ra\infty$ in \eqref{e:prod} in the sense of
 distribution to obtain
\begin{equation}
\label{e:prod.wk}
 \frac{d}{dt}\lra{u_1}{u_2} = \lra{u_{1,t}}{u_2} + \lra{u_{2,t}}{u_1},
\end{equation}
 in the sense of distribution. Notice that
\[ \lra{u_{1,t}}{u_2} \in L^2(0, T),\quad \lra{u_{2,t}}{u_1}\in L^1(0,T). \]
 Therefore, $\lra{u_1}{u_2}\in W^{1,1}(0,T)$. Thus, it is absolutely continuous
 in $t$ and for any $t_0\in[0,T)$ and $t\in(t_0, T]$,
\begin{equation}
\label{e:prod.int}
 \lra{u_1(t)}{u_2(t)} = \lra{u_1(t_0)}{u_2(t_0)}
 + \int_{t_0}^t\lf(\lra{u_{1,t}(s)}{u_2(s)}+\lra{u_{2,t}(s)}{u_1(s)}\rt) ds.
\end{equation}
 By the regularity of $(u_i,v_i,\tht_i)$, $i=1,2$, we have as well
\begin{equation}
\label{e:ut}
 u_{i,t} = -A_1u_i - (u_i,w_i)\cdot\nbl u_i +v_i -\int_z^0\tht_{i,x},\quad
 i =1,2.
\end{equation}
 From \eqref{e:prod.int}-\eqref{e:ut}, we obtain, for any $t_0\in[0, T)$ and
 $t\in(t_0, T]$,
\begin{equation}
\label{e:u1.u2}
\begin{split}
  \lra{u_1(t)}{u_2(t)} +& 2 \int_{t_0}^t a_1(u_1(s),u_2(s)) ds\\
 +& \int_{t_0}^t
 \lf(\lra{(u_1,w_1)\cdot\nbl u_1-v_1+\int_z^0\tht_{1,x}}{u_2}\rt)ds\\
 +& \int_{t_0}^t
 \lf(\lra{(u_2,w_2)\cdot\nbl u_2-v_2+\int_z^0\tht_{1,x}}{u_1}\rt)ds\\
 =&\lra{u_1(t_0)}{u_2(t_0)} .
\end{split}
\end{equation}
 Similarly, we can also obtain, for any $t_0\in[0, T)$ and $t\in(t_0, T]$,
\begin{equation}
\label{e:v1.v2}
\begin{split}
 \lra{v_1(t)}{v_2(t)}
 +& 2 \int_{t_0}^t a_2(v_1(s),v_2(s)) ds\\
 +& \int_{t_0}^t
 \lf(\lra{(u_1,w_1)\cdot\nbl v_1+u_1}{v_2}\rt)ds\\
 +& \int_{t_0}^t
 \lf(\lra{(u_2,w_2)\cdot\nbl v_2+u_2}{v_1}\rt)ds\\
 =&\lra{v_1(t_0)}{v_2(t_0)},
\end{split}
\end{equation}
 and for any $t_0\in[0, T)$ and $t\in(t_0, T]$,
\begin{equation}
\label{e:t1.t2}
\begin{split}
 \lra{\tht_1(t)}{\tht_2(t)}
 +& 2 \int_{t_0}^t a_3(\tht_1(s),\tht_2(s)) ds\\
 +& \int_{t_0}^t \lra{(u_1,w_1)\cdot\nbl\tht_1}{\tht_2} ds\\
 +& \int_{t_0}^t \lra{(u_2,w_2)\cdot\nbl\tht_2}{\tht_1} ds\\
 =&\lra{\tht_1(t_0)}{\tht_2(t_0)} +\int_{t_0}^t \lra{Q}{\tht_1+\tht_2}ds.
\end{split}
\end{equation}
 By Defintion~\ref{d:soln}, $(u_2,v_2,\tht_2)$ satifies the energy inequalities
 \eqref{e:ei.u}-\eqref{e:ei.t} for {\em almost every} $t\in (0,T]$.
 By Theorem~\ref{t:ee.ss}, $(u_1,v_1,\tht_1)$ satifies the energy equalities
 \eqref{e:ee.u}-\eqref{e:ee.t} for any $t_0\in[0, T)$ and $t\in(t_0, T]$.
 Combining these with \eqref{e:u1.u2}-\eqref{e:t1.t2}, we obtain, for
 {\em almost every} $t\in(0, T]$,
\begin{equation}
\label{e:ie.tu}
 \|\tu(t)\|^2 + 2\int_0^t \|\tu\|_{V_1}^2ds
 \ls -2\int_0^t\lra{\tu u_{1,x} + \tw u_{1,z}
 +\int_z^0\tth_xd\zt + \tv}{\tu} ds,
\end{equation}
\begin{equation}
\label{e:ie.tv}
 \|\tv(t)\|^2 +2 \int_0^t \|\tv\|_{V_2}^2ds
 \ls -2 \int_0^t \lra{\tu v_{1,x} +\tw v_{1,z} -\tu}{\tv} ds,
\end{equation}
\begin{equation}
\label{e:ie.tt}
 \|\tth(t)\|^2 +2\int_0^t \|\tth\|_{V_3}^2 ds
 \ls -2\int_0^t \lra{\tu \tht_{1,x} + \tw \tht_{1,z}}{\tth} ds.
\end{equation}
 Notice that some cancellations are used in the derivation of the above
 inequalities, which can be justified using Lemma~\ref{l:tri}. We omit
 justification of these cancellations here, since we have done similar
 justifications before.
 Now, we estimate the terms on the right-hand side of equations
 \eqref{e:ie.tu}-\eqref{e:ie.tt}. By Lemma~\ref{l:tri}, we have
\begin{equation}
\label{e:nlt}
\begin{split}
 \lf|\lra{\tu u_{1,x}}{\tu}\rt|
\lsp& \|u_{1,x}\|\|\tu\|\|\tu\|_{V_1},\\
 \lf|\lra{\tw u_{1,z}}{\tu}\rt|
\lsp& \|u_{1,z}\|\|\tu\|^\hf\|\tu\|_{V_1}^{\frac{3}{2}},\\
 \lf|\lra{\tu v_{1,x}}{\tv}\rt|
\lsp& \|v_{1,x}\|\|\tu\|^\hf\|\tu\|_{V_1}^\hf\|\tv\|^\hf\|\tv\|_{V_2}^\hf,\\
 \lf|\lra{\tw v_{1,z}}{\tv}\rt|
\lsp& \|v_{1,z}\| \|\tu\|_{V_1}\|\tv\|^\hf\|\tv\|_{V_2}^\hf,\\
 \lf|\lra{\tu \tht_{1,x}}{\tth}\rt|
\lsp& \|\tht_{1,x}\|\|\tu\|^\hf\|\tu\|_{V_1}^\hf\|\tth\|^\hf\|\tth\|_{V_3}^\hf,\\
 \lf|\lra{\tw \tht_{1,z}}{\tth}\rt|
\lsp& \|\tht_{1,z}\| \|\tu\|_{V_1}\|\tth\|^\hf\|\tth\|_{V_3}^\hf.
\end{split}
\end{equation}
 Summing up \eqref{e:ie.tu}-\eqref{e:ie.tt}, applying the estimates in
 \eqref{e:nlt} and using Cauchy-Schwartz inequality, we get for {\em almost
 every} $t\in(0, T]$,
\begin{equation*}
\begin{split}
 \|(\tu(t),\tv(t),\tth(t))\|_{H}^2 &+\int_0^t\|(\tu(t),\tv(t),\tth(t))\|_{V}^2\\
 &\lsp \int_0^t (1+\|(u_1,v_1, \tht_1)\|_V^4)\|(\tu,\tv,\tth)\|_{H}^2 ds.
\end{split}
\end{equation*}
 Applying a generalized version of Gronwall lemma to the above inequality yields
\[ \|(\tu(t), \tv(t),\tth(t))\|_{H} = 0, \quad\ \text{for a.e.}\ 
 t\in [0,T]. \]
 Therefore, $(u_1,v_1,\tht_1) = (u_2,v_2,\tht_2)$, for almost every $t\in[0,T]$.

\end{prf}

\begin{rmrk}
\label{r:lc}
{\em We have in fact proved, for general $(\tu(0), \tv(0),\tth(0))\in H$ and
 $(u_1(0),v_1(0),\tht_1(0))\in V$, the following Lipschitz continuity property
 for every\footnote{See \eqref{e:wscb.H}.} $t\in[0,T]$, 
\begin{equation}
\begin{split}
 &\|(\tu(t), \tv(t),\tth(t))\|_H^2\\
 \lsp& \|(\tu(0), \tv(0),\tth(0))\|_{H}^2
 \exp\lf\{\int_0^t (1+\|(u_1,v_1, \tht_1)\|_V^4) ds\rt\}.
\end{split}
\end{equation}
}
\done
\end{rmrk}

 The following theorem gives a much deeper description of a weak solution
 than Definition~\ref{d:soln} and Theorem~\ref{t:ws} combined.
\begin{thrm}
\label{t:ee.ws}
 Let $(0,T)$ be the largest interval of existence for a weak solution
 $(u,v,\tht)$ of the problem \eqref{e:u}-\eqref{e:ic}
 with $(u_0,v_0,\tht_0)\in H$. Then, $T=\infty$. Moreover,
\begin{equation}
\label{e:wsc.V}
  (u,v,\tht) \in C((0,\infty), V),
\end{equation}
\begin{equation}
\label{e:h2.ae}
 (u(t),v(t),\tht(t))\in D(A), \quad \text{for a.e.}\ t>0,
\end{equation}
 and, for any $t_0\in[0,\infty)$, $t\in(t_0, \infty)$, the energy equalities
 \eqref{e:ee.u}-\eqref{e:ee.t} are valid. Moreover,
\begin{equation}
\label{e:lsc.H.2}
 \lim_{t\ra 0^+}\|(u(t),v(t),\tht(t))-(u_0,v_0,\tht_0)\|_H =0.
\end{equation}
 Therefore,
\begin{equation}
\label{e:wscb.H} (u,v,\tht) \in C_B([0,\infty), H).
\end{equation}
\end{thrm}

\begin{prf}
 Let $(0,T)$ be the largest interval of existence for a weak solution
 $(u,v,\tht)$ of \eqref{e:u}-\eqref{e:ic}. Since $(u,v,\tht)\in L^2(0,T;V)$,
 we have
\[ (u(t), v(t),\tht(t))\in V, \quad \text{for a.e.}\ t\in(0, T]. \]
 Choose $\tau\in (0, T)$ such that $(u(\tau),v(\tau),\tht(\tau))\in V$ and that
 \eqref{e:ei.u}-\eqref{e:ei.t} are satisfied with $t_0=\tau$. Then, by
 Remark~\ref{r:ss}, there is a strong solution $(u_1,v_1,\tht_1)$ of
 \eqref{e:u}-\eqref{e:ic} on $[\tau, \infty)$ such that
\[ (u_1(\tau),v_1(\tau),\tht_1(\tau)) = (u(\tau),v(\tau),\tht(\tau)). \]
 By Definition~\ref{d:soln}, $(u,v,\tht)$ is a weak solution on $[\tau,T]$.
 Therefore, by Theorem~\ref{t:ws.uniq}, $(u_1,v_1,\tht_1)=(u,v,\tht)$ on
 $[\tau, T]$. By Theorem~\ref{t:ee.ss}, $(u_1,v_1,\tht_1)$ is also a weak
 solution on $[\tau,\infty)$. By Theorem~\ref{t:ws.uniq} again,
 $(u_1,v_1,\tht_1)$ is the unique weak solution on $[\tau,\infty)$. This proves
 $T=\infty$.
 
 Notice that the above $\tau$ can be chosen arbitrarily small. Therefore,
 \eqref{e:wsc.V} follows from continuity property \eqref{e:ss.c} for a strong
 solution, \eqref{e:h2.ae} follows from the definition of a strong solution,
 and by Theorem~\ref{t:ee.ss}, for any $t_0>0$ and all $t\in(t_0,\infty)$, the
 energy equalities \eqref{e:ee.u}-\eqref{e:ee.t} are satisfied.
 
 Next, we prove validity of \eqref{e:ee.u}-\eqref{e:ee.t} for $t_0=0$ and all
 $t>0$. By Theorem~\ref{t:ws}, there exists a set $E\subset(0,\infty)$ such
 that \eqref{e:lsc.H} is valid. So, we can choose a sequence
\[ \{t_n\}_{n=1}^\infty \subset (0,\infty) \setminus E, \]
 which is monotonically decreasing to $0$ as $n\ra\infty$ and 
\begin{equation}
\label{e:lsc}
 \lim_{n\ra\infty}\|(u(t_n),v(t_n),\tht(t_n))-(u_0,v_0,\tht_0)\|_H =0.
\end{equation}
 Since $t_n>0$ for every $n\gs1$, we have just proved, for any $t>t_n$,
\begin{equation}
\label{e:ee.u.tn}
\|u(t)\|^2+2\int_{t_n}^t\lf(\|u(s)\|_{V_1}^2
 + \lra{-v(s)+\int_z^0\tht_x(s)}{u(s)} \rt) ds =\|u(t_n)\|^2,
\end{equation}
\begin{equation}
\label{e:ee.v.tn}
\|v(t)\|^2+2\int_{t_n}^t \lf(\|v(s)\|_{V_2}^2 +\lra{u(s)}{v(s)}\rt)ds
 =\|v(t_n)\|^2,
\end{equation}
\begin{equation}
\label{e:ee.t.tn}
\|\tht(t)\|^2 + \int_{t_n}^t\lf(\|\tht(s)\|_{V_3}^2-\lra{Q}{\tht(s)}\rt) ds =
 \|\tht(t_n)\|^2.
\end{equation}
 Now, take the limit $n\ra\infty$ in \eqref{e:ee.u.tn}-\eqref{e:ee.t.tn} and
 using the continuity property \eqref{e:lsc} and the fact that
 $(u,v,\tht)\in L^2(0,\infty;V)$, we have, for any $t>0$,
\begin{equation}
\label{e:ee.u.0}
\|u(t)\|^2+2\int_0^t \lf(\|u(s)\|_{V_1}^2
 + \lra{-v(s)+\int_z^0\tht_x(s)}{u(s)} \rt) ds =\|u_0\|^2,
\end{equation}
\begin{equation}
\label{e:ee.v.0}
\|v(t)\|^2+2\int_0^t \lf(\|v(s)\|_{V_2}^2 +\lra{u(s)}{v(s)}\rt)ds
 =\|v_0\|^2,
\end{equation}
\begin{equation}
\label{e:ee.t.0}
\|\tht(t)\|^2 + \int_0^t\lf(\|\tht(s)\|_{V_3}^2-\lra{Q}{\tht(s)}\rt) ds =
 \|\tht_0\|^2.
\end{equation}
 These are \eqref{e:ee.u}-\eqref{e:ee.t} for $t_0=0$ and all $t>0$.

 Moreover, by \eqref{e:ee.u.0} and that $(u,v,\tht)\in L^2(0,\infty;V)$, we
 obtain
\[ \limsup_{t\ra 0^+} \|u(t)\|^2 \ls \|u_0\|^2. \]
 By \eqref{e:c.w}, we also have
\[ \|u_0\|^2 \ls \liminf_{t\ra 0^+} \|u(t)\|^2. \]
 Thus,
\[ \lim_{t\ra 0^+} \|u(t)\|^2 = \|u_0\|^2. \]
 Using \eqref{e:c.w} again, we obtain
\[ \lim_{t\ra 0^+} \|u(t)-u_0\|^2=0. \]
 Similarly, we have
\[ \lim_{t\ra 0^+} \|v(t)-v_0\|^2=\lim_{t\ra 0^+} \|\tht(t)-\tht_0\|^2=0. \]
 This proves \eqref{e:lsc.H.2}.

 Finally, by definition, as a weak solution,
\[ (u,v,\tht)\in L^\infty(0, T_0;H),\quad \forall  T_0>0. \]
 As a strong solution, the uniform boundedness in $V$ is valid:
\[ (u,v,\tht)\in L^\infty(T_0, \infty; V). \]
 Therefore,
\[ (u,v,\tht) \in L^\infty(0,\infty; H). \]
 Then, \eqref{e:wscb.H} follows form the above uniform boundedness in $H$,
 \eqref{e:wsc.V} and \eqref{e:lsc.H.2}.

\end{prf}

\section{A Sufficient Condition for Uniqueness}
\label{s:sc}

 In this section, we present a new sufficient condition for uniqueness of weak
 solutions of 2D viscous PE \eqref{e:u}-\eqref{e:ic}. First, we mention the
 following result for a sufficient condition for uniqueness of weak solutions
 of \eqref{e:u}-\eqref{e:ic}:
\begin{prop}
\label{p:sc}
 Let $(u_i, v_i, \tht_i)$, for $i=1,2$, be weak solutions of
 \eqref{e:u}-\eqref{e:ic}. Suppose $(u_0, v_0, \tht_0)\in H$ and for some $T>0$,
\begin{equation*}
 (u_{1,z}, v_{1,z}, \tht_{1,z}) \in L^4(0, T; [L^2(D)]^3).
\end{equation*}
 Then, $(u_1,v_1,\tht_1)\equiv(u_2,v_2,\tht_2)$ for $t\in [0, T]$.
\end{prop}

 Proposition~\ref{p:sc} was proved in \cite{g.m.r:01} for the 2D hydrostatic
 Navier-Stokes equations, where $v$, $\tht$ are neglected. In \cite{g.m.r:01},
 the definition of weak solution is somewhat different from
 Definition~\ref{d:soln} and the boundary condition is also somewhat different.

 As the first main result of this section, the following Theorem~\ref{t:sc}
 generalizes Proposition~\ref{p:sc} and allows one to find new classes of weak
 solutions of the system of \eqref{e:u}-\eqref{e:ic}, within which the weak
 solutions are unique. Especially, it is crucial for proving our main uniqueness
 result in Section~\ref{s:un}.

\begin{thrm}
\label{t:sc}
 Let $(u_i, v_i, \tht_i)$, for $i=1,2$, be weak solutions of
 \eqref{e:u}-\eqref{e:ic}. Suppose $(u_0, v_0, \tht_0)\in H$ and for some $T>0$,
\begin{equation}
\label{e:sc}
(u_{1,z},v_{1,z},\tht_{1,z})
 \in \lf[L^4(0,T;L^2(D))\cup L^2(0,T; L_x^\infty(L_z^2))\rt]^3.
\end{equation}
 Then, $(u_1,v_1,\tht_1)\equiv(u_2,v_2,\tht_2)$ for all $t\in [0, T]$.
\end{thrm}

\begin{prf}

 We will prove Lipschitz continuity of the weak solutions with respect to
 initial data in $L^2$, assuming $(u_1,v_1,\tht_1)$ satisfies the regularity
 condition \eqref{e:sc}.

 Denote:
\[ \tld{u} = u_1-u_2, \quad \tld{v} = v_1-v_2, \quad
  \tld{\tht} = \tht_1-\tht_2,\]
\[ \tld{u}_0 = u_{1,0}-u_{2,0}, \quad \tld{v}_0 = v_{1,0}-v_{2,0}, \quad
  \tld{\tht}_0 = \tht_{1,0}-\tht_{2,0},\]
 and
\[ \tld{w}=w_1-w_2, \quad \tld{q} = q_1 -q_2.\]

 Let $t_0\in(0,T)$. Then $(u_1,v_1,\tht_1)$ and $(u_2,v_2,\tht_2)$ are both
 strong solutions on $[t_0, T]$. Therefore, we can follow the proof of
 Theorem~\ref{t:ws.uniq} and apply it to $(u_1,v_1,\tht_1)$ and
 $(u_2,v_2,\tht_2)$ on $[t_0, T]$ to obtain, for {\em every} $t\in[t_0, T]$,
\begin{equation}
\label{e:ie.tu.t0}
\begin{split}
 \|\tu(t)\|^2 &+ 2\int_{t_0}^t \|\tu\|_{V_1}^2ds\\
 =& \|\tu(t_0)\|^2 -2\int_{t_0}^t\lra{\tu u_{1,x} + \tw u_{1,z}
 +\int_z^0\tth_xd\zt + \tv}{\tu} ds,
\end{split}
\end{equation}
\begin{equation}
\label{e:ie.tv.t0}
 \|\tv(t)\|^2 +2 \int_{t_0}^t \|\tv\|_{V_2}^2ds
 = \|\tv(t_0)\|^2 -2 \int_{t_0}^t \lra{\tu v_{1,x} +\tw v_{1,z} -\tu}{\tv} ds,
\end{equation}
\begin{equation}
\label{e:ie.tt.t0}
 \|\tth(t)\|^2 +2\int_{t_0}^t \|\tth\|_{V_3}^2 ds
 = \|\tth(t_0)\|^2 -2\int_{t_0}^t \lra{\tu \tht_{1,x} + \tw \tht_{1,z}}{\tth} ds.
\end{equation}
 Notice that {\em equalities} \eqref{e:ie.tu.t0}-\eqref{e:ie.tt.t0} are
 obtained similarly to {\em inequalities} \eqref{e:ie.tu}-\eqref{e:ie.tt}.
 Since $(u_1,v_1,\tht_1)$ and $(u_2,v_2,\tht_2)$ are both strong solutions on
 $[t_0, T]$, we indeed have above {\em equalities} and they are valid for
 {\em every} $t\in[t_0, T]$.

 By Theorem~\ref{t:ee.ws}, $(u_i,v_i,\tht_i)\in C([t_0,T], V)$ for $i=1,2$.
 Thus, \eqref{e:ie.tu.t0}-\eqref{e:ie.tt.t0} along with Lemma~\ref{l:tri}
 imply $\|\tu\|, \|\tv\|, \|\tth\|\in C^1([t_0, T])$.
 Therefore, we have in {\em classic sense for every} $t\in[t_0,T]$ that
\begin{equation}
\label{e:tu.2}
\begin{split}
 \hf\dt{\|\tu\|^2} +\|\tu\|_{V_1}^2 = -\lra{\tu u_{1,x} + \tw u_{1,z}
 +\int_z^0\tth_xd\zt + \tv}{\tu},
\end{split}
\end{equation}
\begin{equation}
\label{e:tv.2}
\begin{split}
 \hf\dt{\|\tv\|^2} +\|\tv\|_{V_2}^2
 = -\lra{v_{1,x}\tu +\tw v_{1,z} -\tu}{\tv},
\end{split}
\end{equation}
\begin{equation}
\label{e:tt.2}
\begin{split}
 \hf\dt{\|\tth\|^2} +\|\tth\|_{V_3}^2
 = -\lra{\tu \tht_{1,x} + \tw \tht_{1,z}}{\tth}.
\end{split}
\end{equation}
 Summing up \eqref{e:tu.2}-\eqref{e:tt.2} yields, for every $t\in[t_0,T]$,
\begin{equation}
\label{e:tld.2}
\begin{split}
\hf\frac{d}{dt} &\lf(\|\tu\|^2+\|\tv\|^2+\|\tth\|^2\rt)
  +\|\tu\|_{V_1}^2+\|\tv\|_{V_2}^2+\|\tth\|_{V_3}^2\\
 =&-\lra{\tu u_{1,x} + \tw u_{1,z}
 +\int_z^0\tth_xd\zt}{\tu}\\
 &-\lra{v_{1,x}\tu +\tw v_{1,z}}{\tv}
 -\lra{\tu \tht_{1,x} + \tw \tht_{1,z}}{\tth}
\end{split}
\end{equation}
 Now, we estimates all the terms on the right side of \eqref{e:tld.2}.
 The bilinear term is easily estimated by Cauchy-Schwartz inequality:
\begin{equation}
\label{e:t.u}
 \lf|\lra{\int_z^0\tth_xd\zt}{\tu}\rt| \ls \|\tth_x\|\|\tu\|
\ls C_\veps\|\tu\|^2+ \veps\|\tth_x\|^2.
\end{equation}
 Next, by Agmon's inequality, we have
\begin{equation*}
\begin{split}
\lf|\lra{\tu u_{1,x}}{\tu}\rt|
 \ls& \int_0^1 \|\tu\|_{L^\infty_z}\lf(\int_{-h}^0 |u_{1,x} \tu| dz\rt) dx \\
 \ls& \int_0^1 \|\tu\|_{L^2_z}^\hf \|\tu_z\|_{L^2_z}^\hf \|u_{1,x}\|_{L^2_z}
  \|\tu\|_{L^2_z} dx \\
 \ls& \|\tu\|_{L^\infty_x(L^2_z)}
 \int_0^1 \|\tu\|_{L^2_z}^\hf \|\tu_z\|_{L^2_z}^\hf \|u_{1,x}\|_{L^2_z} dx\\
\ls& \|\tu\|_{L^\infty_x(L^2_z)}\|\tu\|^\hf \|\tu_z\|^\hf \|u_{1,x}\|,
\end{split}
\end{equation*}
 where H\"older's inequality is applied in the last step. Then, by Minkowski's
 inequality:
\[ \|\tu\|_{L^\infty_x(L^2_z)} \ls \|\tu\|_{L^2_z(L^\infty_x)}, \]
 we get
\begin{equation}
\label{e:u.n0}
\begin{split}
\lf|\lra{\tu u_{1,x}}{\tu}\rt|
 \ls& \|\tu\|_{L^2_z(L^\infty_x)} \|\tu\|^\hf \|\tu_z\|^\hf \|u_{1,x}\|\\ 
 \ls& \text{ (Agmon's inequality) }\\
 \ls& \|u_{1,x}\| \|\tu\|^\hf\|\tu_x\|^\hf\|\tu\|^\hf\|\tu_z\|^\hf \\
 =& \|u_{1,x}\| \|\tu\|\|\tu_x\|^\hf\|\tu_z\|^\hf\\
 \ls& C_\veps\|u_{1,x}\|^2\|\tu\|^2 + \veps \|\nbl\tu\|^2.
\end{split}
\end{equation}
 It follows similarly that
\begin{equation}
\label{e:v.n0}
\begin{split}
\lf|\lra{\tu v_{1,x}}{\tv}\rt|
\ls& \|v_{1,x}\| \|\tu\|^\hf\|\tu_x\|^\hf\|\tv\|^\hf \|\tv_z\|^\hf\\
\ls& C_\veps\|v_{1,x}\|^2 (\|\tu\|^2+\|\tv\|^2) + \veps(\|\tu_x\|^2+\|\tv_z\|^2),
\end{split}
\end{equation}
 and
\begin{equation}
\label{e:t.n0}
\begin{split}
\lf|\lra{\tu \tht_{1,x}}{\tth}\rt|
\ls& \|\tht_{1,x}\| \|\tu\|^\hf\|\tu_x\|^\hf\|\tth\|^\hf \|\tth_z\|^\hf\\
\ls& C_\veps\|\tht_{1,x}\|^2 (\|\tu\|^2+\|\tth\|^2)
 + \veps(\|\tu_x\|^2+\|\tth_z\|^2).
\end{split}
\end{equation}
 Finally, similar to above estimates, we have
\begin{equation*}
\begin{split}
\lf|\lra{\tw u_{1,z}}{\tu}\rt|
 \ls& \int_0^1 \|\tw\|_{L^\infty_z} \lf(\int_{-h}^0 |u_{1,z} \tu| dz\rt) dx\\
\ls&\int_0^1 \|\tw\|_{L^2_z}^\hf\|\tw_z\|_{L^2_z}^\hf
 \|u_{1,z}\|_{L^2_z}\|\tu\|_{L^2_z} dx\\
\lsp&\int_0^1 \|\tu_x\|_{L^2_z}\|u_{1,z}\|_{L^2_z}\|\tu\|_{L^2_z} dx.
\end{split}
\end{equation*}
 The right-hand side of the above inequality can be further estimated in two
 different ways:
\begin{equation}
\label{e:u.n1}
\begin{split}
\lf|\lra{\tw u_{1,z}}{\tu}\rt|
\ls& C\|\tu\|_{L^\infty_x(L^2_z)}\int_0^1 \|\tu_x\|_{L^2_z}\|u_{1,z}\|_{L^2_z} dx\\
\ls& C\|\tu\|_{L^2_z(L^\infty_x)} \|u_{1,z}\| \|\tu_x\|\\
\ls& C\|\tu\|^\hf \|\tu_x\|^\hf\|u_{1,z}\| \|\tu_x\|\\
\ls& C_\veps\|u_{1,z}\|^4\|\tu\|^2 + \veps \|\tu_x\|^2
\end{split}
\end{equation}
 and
\begin{equation}
\label{e:u.n2}
\begin{split}
\lf|\lra{\tw u_{1,z}}{\tu}\rt|
 \ls& C\|u_{1,z}\|_{L^\infty_x(L^2_z)}
 \int_0^1 \|\tu_x\|_{L^2_z}\|\tu\|_{L^2_z} dx\\
\ls& C\|u_{1,z}\|_{L^\infty_x(L^2_z)} \|\tu_x\| \|\tu\|\\
 \ls& C_\veps\|u_{1,z}\|_{L^\infty_x(L^2_z)}^2 \|\tu\|^2 + \veps \|\tu_x\|^2.
\end{split}
\end{equation}
 In the above estimates \eqref{e:u.n1} and \eqref{e:u.n2}, we have used ideas
 similar to those used in \eqref{e:u.n0}. It now follows, similar to
 \eqref{e:u.n1} and \eqref{e:u.n2}, that
\begin{equation}
\label{e:v.n1}
\begin{split}
\lf|\lra{\tw v_{1,z}}{\tv}\rt|
\ls& C\|\tv\|_{L^\infty_x(L^2_z)}\int_0^1 \|\tu_x\|_{L^2_z}\|v_{1,z}\|_{L^2_z} dx\\
\ls& C\|\tv\|_{L^2_z(L^\infty_x)} \|v_{1,z}\| \|\tu_x\|\\
\ls& C\|\tv\|^\hf \|\tv_x\|^\hf\|v_{1,z}\| \|\tu_x\|\\
\ls& C_\veps\|v_{1,z}\|^4\|\tv\|^2 + \veps (\|\tu_x\|^2+\|\tv_x\|^2),
\end{split}
\end{equation}
\begin{equation}
\label{e:v.n2}
\begin{split}
\lf|\lra{\tw v_{1,z}}{\tv}\rt|
 \ls& C\|v_{1,z}\|_{L^\infty_x(L^2_z)}
 \int_0^1 \|\tu_x\|_{L^2_z}\|\tv\|_{L^2_z} dx\\
\ls& C\|v_{1,z}\|_{L^\infty_x(L^2_z)} \|\tu_x\| \|\tv\|\\
 \ls& C_\veps\|v_{1,z}\|_{L^\infty_x(L^2_z)}^2 \|\tv\|^2 + \veps \|\tu_x\|^2,
\end{split}
\end{equation}
 and
\begin{equation}
\label{e:t.n1}
\begin{split}
\lf|\lra{\tw \tht_{1,z}}{\tth}\rt|
 \ls& C_\veps\|\tth_{1,z}\|^4\|\tu\|^2 + \veps (\|\tu_x\|^2+\|\tth_x\|^2),
\end{split}
\end{equation}
\begin{equation}
\label{e:t.n2}
\begin{split}
\lf|\lra{\tw \tht_{1,z}}{\tth}\rt|
 \ls& C_\veps\|\tht_{1,z}\|_{L^\infty_x(L^2_z)}^2 \|\tth\|^2 + \veps \|\tu_x\|^2.
\end{split}
\end{equation}
 Plug the estimates \eqref{e:t.u}-\eqref{e:t.n0}, \eqref{e:u.n1} or
 \eqref{e:u.n2}, \eqref{e:v.n1} or \eqref{e:v.n2}, and \eqref{e:t.n1} or
 \eqref{e:t.n2} into \eqref{e:tld.2} and choose sufficiently small $\veps>0$ to
 absorb the dissipation terms. Then, assume that $(u_1,v_1,\tht_1)$ satisfies
 \eqref{e:sc} and apply Gronwall lemma to finish the proof. Notice that global
 regularity result of weak solutions of \eqref{e:u}-\eqref{e:ic} is also used in
 justifying applicability of Gronwall lemma. As an example to demonstrate the
 details, we now finish the proof for a special case of \eqref{e:sc} when
 $(u_1,v_1,\tht_1)$ satisfies:
\begin{equation}
\label{e:asmp.ex}
 (u_{1,z},v_{1,z},\tht_{1,z})\in \lf[L^2(0, T; L_x^\infty(L_z^2))\rt]^3.
\end{equation}
 Plug \eqref{e:t.u}-\eqref{e:t.n0}, \eqref{e:u.n2},
 \eqref{e:v.n2} and \eqref{e:t.n2} into \eqref{e:tld.2} and choose sufficiently
 small $\veps>0$, we obtain, for $t\in [t_0, T]$,
\begin{equation}
\begin{split}
 \frac{d}{dt}& \|(\tu,\tv,\tth)\|_H^2
  +\|(\tu,\tv,\tth)\|_V^2\\
  \lsp & \lf[1 + \|(u_1, v_1,\tht_1)_x\|_H^2
 + \|(u_1,v_1,\tht_1)_z\|_{(L_x^\infty(L_z^2))^3)}^2\rt]
 \|(\tu,\tv,\tth)\|_H^2.
\end{split}
\end{equation}
 Noticing \eqref{e:asmp.ex}, we can use Gronwall inequality to obtain, for
 $t\in[t_0, T]$,
\begin{equation*}
\begin{split}
 &\|(\tu(t),\tv(t),\tth(t))\|_H^2\\
\lsp& \|(\tu(t_0),\tv(t_0),\tth(t_0))\|_H^2\\
 &\x \exp\lf\{ \int_0^t \lf[1 + \|(u_1, v_1,\tht_1)_x\|_H^2
 + \|(u_1,v_1,\tht_1)_z\|_{(L_x^\infty(L_z^2))^3)}^2\rt] ds \rt\}.
\end{split}
\end{equation*}
 Now, take the limit $t_0\ra 0^+$ and use Theorem~\ref{t:ee.ws}, we get
\begin{equation*}
\begin{split}
 &\|(\tu(t),\tv(t),\tth(t))\|_H^2\\
\lsp& \|(\tu_0,\tv_0,\tth_0)\|_H^2\\
 &\x \exp\lf\{ \int_0^t \lf[1 + \|(u_1, v_1,\tht_1)_x\|_H^2
 + \|(u_1,v_1,\tht_1)_z\|_{(L_x^\infty(L_z^2))^3)}^2\rt] ds \rt\}.
\end{split}
\end{equation*}
 Since the above inequality is independent of $t_0$ and $t_0$ can be chosen
 arbitrarily small, it is valid for all $t\in(0, T]$.  This proves that
 $(u_1,v_1,\tht_1)\equiv(u_2,v_2,\tht_2)$, if $(\tu_0,\tv_0,\tth_0)=(0,0,0)$.
 The other cases covered by \eqref{e:sc} can be similarly proved. 

\end{prf}

\section{Global Existence}
\label{s:ex}

 In this section, we prove global in time {\em uniform boundedness} of the
 norms of some partial derivatives of the weak solutions. These results are
 also important for proving our uniqueness result in Section~\ref{s:un}.

 Let us mention first that it is not immediately clear whether or not the
 global uniform $L^2_xH^\hf_z$ boundedness for the solutions of the 2D
 hydrostatic Navier-Stokes equations as obtained in \cite{b.k.l:04} can be
 extended to the problem of \eqref{e:u}-\eqref{e:ic}. This is due to the fact
 that the boundary conditions \eqref{e:bc.u} and \eqref{e:bc.v} for $(u,v)$
 are different from the boundary condition \eqref{e:bc.t} for $\tht$ and the
 possibility that $\alf_1$ and $\alf_2$ may be different. This problem will be
 studied elsewhere.

 We begin with a theorem for global in time uniform boundedness of
 $(u_z,v_z,\tht_z)$ in $[L^2(D)]^3$ for \eqref{e:u}-\eqref{e:ic}.
 Recall that global existence and uniqueness of $z$-weak solutions were proved
 in \cite{p:07} for 2D viscous PE in the case of periodic boundary conditions;
 and in \cite{j:17} for 3D viscous PE in case of Neumann boundary condition
 for horizontal velocity at bottom of the physical domain.
 However, these analyses do not apply to the system \eqref{e:u}-\eqref{e:ic}
 due to different boundary conditions. A possible approach might be a proper
 modification of that of \cite{b.k.l:04} in obtaining boundedness for $z$-weak
 solutions of the simplified 2D hydrostatic Navier-Stokes equations, where $v$
 and $\tht$ were omitted. Nevertheless, new issues will come up again due to
 boundary conditions.
 Instead, we will take advantage of a result of \cite{p.t.z:09} directly in our
 proof of the following Theorem~\ref{t:uz}.

 \begin{thrm}
\label{t:uz}
 Suppose $Q\in L^2(D)$, $(u_0, v_0,\tht_0)\in H$ and $(u,v,\tht)$ is a weak
 solution of \eqref{e:u}-\eqref{e:ic}.  The following statements are valid:
\begin{enumerate}
 \item[{\em(a)}] If $\Pz u_0\in L^2(D)$, then there exists a weak solution
 $(u,v,\tht)$ of \eqref{e:u}-\eqref{e:ic}, such that
 \[ u_z\in L^\infty(0,\infty; L^2(D)) \cap L^2(0,\infty; H^1(D)). \]
 Moreover, there exists a bounded absorbing set for $u_z$ in $L^2(D)$.
 \item[{\em(b)}] If $(\Pz u_0,\Pz v_0)\in (L^2(D))^2$, then there exists a weak
 solution $(u,v,\tht)$ of \eqref{e:u}-\eqref{e:ic}, such that
\[ (u_z,v_z)\in L^\infty(0,\infty; [L^2(D)]^2)\cap L^2(0,\infty; [H^1(D)]^2).\]
 Moreover, there exists a bounded absorbing set for $(u_z,v_z)$ in $[L^2(D)]^2$.
 \item[{\em(c)}] If $(\Pz u_0,\Pz\tht_0)\in (L^2(D))^2$, then there exists a
 weak solution $(u,v,\tht)$ of \eqref{e:u}-\eqref{e:ic}, such that
\[ (u_z,\tht_z))\in L^\infty(0,\infty; [L^2(D)]^2)
     \cap L^2(0,\infty; [H^1(D)]^2).\]
 Moreover, there exists a bounded absorbing set for $(u_z,\tht_z)$ in
 $[L^2(D)]^2$.
\end{enumerate}
\end{thrm}

\begin{rmk}
 Quite unexpectedly, it seems to be a non-trivial problem whether or not global
 in time uniform boundedness of $\|(u_z, v_x)\|_{[L^2(D)]^2}$ or
 $\|(u_z, \tht_x)\|_{[L^2(D)]^2}$ is still valid when
 $(u_{0,z}, v_{0,x})\in [L^2(D)]^2$ or $(u_{0,z}, \tht_{0,x})\in [L^2(D)]^2$.

\end{rmk}

\begin{prf}

 {\bf Step 1}. Proof of part (a) of Theorem~\ref{t:uz}.\\
 By Theorem~\ref{t:ee.ws}, we can choose a monotonically decreasing
 sequence
\[ \{ t_n \}_{n=1}^\infty\subset (0, \infty),\text{ such that }\
 \lim_{n\ra\infty} t_n =0, \]
 and
\begin{equation}
\label{e:tn.inf}
 (u,v,\tht)\in C([t_n, \infty), V)\cap L^2(t_n,\infty; D(A)),
 \ \forall\ n\gs 1.
\end{equation}
 Moreover, there exists an absorbing set for $(u,v,\tht)$ in $V$, when
 the time interval $[t_1, \infty)$ is considered.  Therefore,
 {\em what is still needed to be proved is just the following:}
\begin{equation}
\label{e:uz.0.t1}
 u_z\in L^\infty(0, t_1; L^2(D))\cap L^2(0, t_1; H^1(D)).
\end{equation}
 By the estimate of $\|u_z\|$ in \S 3.3 of \cite{p.t.z:09} for a strong
 solution $(u,v,\tht)$ on $[t_n, \infty)$ with initial data
 $(u(t_n),v(t_n),\tht(t_n)\in V$ and by Theorem~\ref{t:ws.uniq}, we have for
 almost every $t\in [t_n,\infty)$,
\begin{equation}
\label{e:uz.ttn.d}
\begin{split}
 \frac{d}{dt} & \lf(\|u_z(t)\|^2 + \alf_1\|u(t)|_{z=0}\|^2\rt)
+ \|\nbl u_z\|^2 + \alf_1\|u_x|_{z=0}\|^2 \\
&\lsp 
 \|\nbl u\|^2+\|v\|^2+\|\tht_x\|^2.
\end{split}
\end{equation}
 Notice that \eqref{e:tn.inf}
 is used in deriving \eqref{e:uz.ttn.d}. Therefore, we have for $t\in[t_n, t_1]$
 with $n>1$,
\begin{equation}
\label{e:uz.ttn}
\begin{split}
 &\|u_z(t)\|^2 + \alf_1\|u(t)|_{z=0}\|^2
+ \int_{t_n}^t \lf(\|\nbl u_z\|^2 + \alf_1\|u_x|_{z=0}\|^2\rt) ds \\
\ls& \|u_z(t_n)\|^2 + \alf_1\|u(t_n)|_{z=0}\|^2 + C \int_{t_n}^t
 \lf[ \|\nbl u\|^2+\|v\|^2+\|\tht_x\|^2\rt] ds.
\end{split}
\end{equation}
 Since $u(t_n)\in V_1$, we have
\begin{equation}
\label{e:u.z=0}
 \lf\|u(t_n)|_{z=0}\rt\| = \lf\|\int_{-h}^0 u_z(t_n) dz\rt\| \lsp
 \|u_z(t_n)\|.
\end{equation}
 Due to the fact that $(u,v,\tht)$ is a weak solution on $(0,\infty)$,
 we also have, for $t\in[t_n, t_1]$,
\begin{equation}
\label{e:uz.irh}
 \int_{t_n}^t (\|\nbl u\|^2+\|v\|^2+\|\tht_x\|^2) ds
 \ls \int_0^{t_1} (\|\nbl u\|^2+\|v\|^2+\|\tht_x\|^2) ds,
\end{equation}
 the upper bound of which depends only on $\|(u_0,v_0,\tht_0)\|_H$, $\|Q\|$
 and $t_1$. Combining \eqref{e:uz.ttn}-\eqref{e:uz.irh}, we have,
 for $t\in[t_n,t_1]$,
\begin{equation}
\label{e:uz.ttn.2}
\begin{split}
 &\|u_z(t)\|^2 + \alf_1\|u(t)|_{z=0}\|^2
+ \int_{t_n}^t \lf(\|\nbl u_z\|^2 + \alf_1\|u_x|_{z=0}\|^2\rt) ds \\
\lsp& \|u_z(t_n)\|^2 + \int_0^{t_1}\lf[ \|\nbl u\|^2+\|v\|^2+\|\tht_x\|^2\rt] ds.
\end{split}
\end{equation}
 Now, choose any $\phi\in \C_0^\infty(D)$. Then, $\phi_z\in V_1$.  Thus, by
 weak continuity of $(u,v,\tht)$ on $[0,\infty)$ (see Theorem~\ref{t:ws}),
 we have
\[ \lim_{n\ra\infty}\lra{u_z(t_n)-\Pz u_0}{\phi}
 = -\lim_{n\ra\infty}\lra{ u(t_n) -u_0}{\phi_z}=0. \]
 Since $\C_0^\infty(D)$ is dense in $L^2(D)$, we have weak convergence:
\[ u_z(t_n) \rau \Pz u_0,\ \text{ in }\ L^2(D). \]
 Therefore, $\{ u_z(t_n) \}_{n=1}^\infty$ is bounded in $L^2(D)$. Now, taking
 the limit $t_n\ra 0$ in \eqref{e:uz.ttn.2}, we have, for all $t\in(0,t_1)$,
\begin{equation}
\label{e:uz.ttn.3}
\begin{split}
 &\|u_z(t)\|^2 + \alf_1\|u(t)|_{z=0}\|^2
+ \int_0^t \lf(\|\nbl u_z\|^2 + \alf_1\|u_x|_{z=0}\|^2\rt) ds \\
\lsp& \sup_{n\gs 1} \|u_z(t_n)\|^2 +
 \int_0^{t_1}\lf[ \|\nbl u\|^2+\|v\|^2+\|\tht_x\|^2\rt] ds.
\end{split}
\end{equation}
 This proves \eqref{e:uz.0.t1} and thus finishes {\bf Step 1}.

 {\bf Step 2}.  Proof of Theorem~\ref{t:uz} part (b).\\
 For simplicity of presentation, in the proof of part (b) of Theorem~\ref{t:uz},
 we will only provide the key estimate of $\|v_z\|$ under the assumption that
 $(u,v,\tht)$ is a {\em strong} solution on $[0, \infty)$. The justification
 that this estimate is sufficient for a rigorous proof of part (b) of
 Theorem~\ref{t:uz} is almost the same as the one we provided in {\bf Step 1}
 for our proof of part (a) of Theorem~\ref{t:uz}. Thus, it is omitted for
 conciseness.

 Taking inner product of \eqref{e:v} with $-v_{zz}$ yields:
\begin{equation}
\label{e:vz.2}
\begin{split}
 \hf\frac{d}{dt}\lf(\|v_z\|^2 + \alf_2\|v(z=0)\|^2\rt) +& \|\nbl v_z\|^2
 +\alf_2\|v_x(z=0)\|^2 \\
 = & \lra{uv_x+wv_z-u}{v_{zz}}.
\end{split}
\end{equation}
 The following computations are used in deriving \eqref{e:vz.2}:
\begin{equation*}
\begin{split}
 -\int v_t v_{zz}
=& -\int_0^1\lf(v_tv_z\Big|_{z=-h}^0 -\int_{-h}^0v_{zt}v_z\rt)dx\\
=& \hf \dt{\|v_z\|^2} + \int_0^1 v_t\alf_2v\Big|_{z=0}\\
=& \hf \dt{\|v_z\|^2} +\frac{\alf_2}{2}\dt{\|v(z=0)\|^2},
\end{split}
\end{equation*}
\begin{equation*}
\begin{split}
 \int_Dv_{zz}v_{xx}
=& \int_{-h}^0\lf(v_{zz}v_x\Big|_{x=0}^1-\int_0^1v_{xzz}v_x dx\rt) dz\\
=&-\int_D v_{xzz}v_x dxdz\\
=&-\int_0^1\lf( v_{xz}v_x\Big|_{z=-h}^0-\int_{-h}^0 v_{xz}^2dz \rt) dx\\
=& \|v_{xz}\|^2+\alf_2\|v_x(z=0)\|^2.
\end{split}
\end{equation*}
 The two tri-linear terms on the right-hand side of \eqref{e:vz.2} will be
 estimated in the following.\\
 First, we have
\begin{equation}
\label{e:wvz.vzz}
\begin{split}
 \int_D wv_z v_{zz} =& \hf \int_D w\Pz(v_z^2)\\
=& \hf\int_0^1 \lf(w v_z^2\Big|_{-h}^0-\int_{-h}^0 w_z v_z^2\rt)dx\\
=& \hf \int_D u_x v_z^2\\
\ls& C\|u_x\|\|v_z\|\|v_{zx}\|^\hf\|v_{zz}\|^\hf\\
\ls& C_\veps\|u_x\|^2\|v_z\|^2 + \veps\|\nbl v_z\|^2.
\end{split}
\end{equation}
 Contrary to the common intuition from experience, the other tri-linear term
 is more complicated to deal with. Integrating by parts and applying boundary
 conditions, one has
\begin{equation*}
\begin{split}
 \int_D uv_x v_{zz} =& \int_0^1
  \lf[ uv_x v_z\Big|_{z=-h}^0-\int_{-h}^0 (u_zv_xv_z+uv_{xz}v_z)dz\rt]dx\\
=&-\alf_2\int_0^1 u v_x v\Big|_{z=0} dx
  - \int_D u_z v_x v_z dxdz - \int_D u v_{xz} v_z dxdz\\
=:& I_0 + I_1 + I_2.
\end{split}
\end{equation*}
 In the following, we estimate $I_0$, $I_1$ and $I_2$ respectively. First,
 we have
\begin{equation*}
\begin{split}
 |I_1|\ls& C\|u_z\| \|v_x\|^\hf\|v_{xz}\|^\hf \|v_z\|^\hf\|v_{zx}\|^\hf\\
  =& C\|u_z\| \|v_x\|^\hf \|v_z\|^\hf \|v_{xz}\|\\
 \ls& C_\veps(\|u_z\|^4 + \|v_x\|^2\|v_z\|^2) + \frac{\veps}{2}\|v_{xz}\|^2,
\end{split}
\end{equation*}
 and
\begin{equation*}
\begin{split}
 |I_2| \ls& C\|v_{xz}\| \|u\|^\hf\|u_z\|^\hf \|v_z\|^\hf \|v_{xz}\|^\hf\\
 =& C\|u\|^\hf\|u_z\|^\hf \|v_z\|^\hf \|v_{xz}\|^{\frac{3}{2}}\\
\ls& C_\veps\|u\|^2\|u_z\|^2\|v_z\|^2 + \frac{\veps}{2}\|v_{xz}\|^2.
\end{split}
\end{equation*}
 Noticing that $u(0,z,t)=0$, we have
\[ u^2(x,0,t) = 2 \int_0^x u(\xi,0,t) u_x(\xi,0,t) d\xi. \]
 Thus,
\[ \|u(z=0)\|_\infty^2 \ls 2 \|u(z=0)\|\|u_x(z=0)\|. \]
 So, we can estimate $I_0$ as following:
\begin{equation}
\label{e:i0}
\begin{split}
 |I_0|
\ls& \alf_2 \|u(z=0)\|_\infty \|v(z=0)\| \|v_x(z=0)\|\\
\ls& \frac{\alf_2}{2}\|u(z=0)\|_\infty^2 \|v(z=0)\|^2
 +\frac{\alf_2}{2}\|v_x(z=0)\|^2\\
\ls& \alf_2\|u(z=0)\| \|u_x(z=0)\| \|v(z=0)\|^2
 + \frac{\alf_2}{2}\|v_x(z=0)\|^2.
\end{split}
\end{equation}
 Due to the boundary condition $u(z=-h)=0$, it holds that,
\[ u(z=0) = \int_{-h}^0 u_z(x,\zt,t)d\zt, \]
 and thus
\begin{equation}
\label{e:be.1}
\begin{split}
 \|u(z=0)\|
=& \lf[ \int_0^1 \lf(\int_{-h}^0 u_z d\zt\rt)^2 dx \rt]^\hf\\
\ls& \lf[ \int_0^1 \lf(\int_{-h}^0 |u_z| d\zt\rt)^2 dx \rt]^\hf\\
\ls& \text{ (Minkowski inequality) } \\
\ls& \int_{-h}^0 \lf(\int_0^1 u_z^2 dx \rt)^\hf dz\\
\ls& h^\hf \|u_z\|.
\end{split}
\end{equation}
 Similar to \eqref{e:be.1}, we also have $\|u_x(z=0)\| \ls h^\hf \|u_{xz}\|$.
 Therefore,
\[ I_0 \ls \alf_2 h \|u_z\| \|u_{xz}\|\|v(z=0)\|^2
 + \frac{\alf_2}{2}\|v_x(z=0)\|^2. \]
 Combining the above estimates of $I_1$, $I_2$ and $I_0$, we have
\begin{equation}
\label{e:uvx.vzz.2}
\begin{split}
 \lf|\int_D uv_xv_{zz} dxdz\rt|
 \ls& C_\veps(\|u\|^2\|u_z\|^2+\|v_x\|^2)\|v_z\|^2\\
 & +\alf_2h\|u_z\|\|u_{xz}\|\|v(z=0)\|^2 + C_\veps \|u_z\|^4\\
 & +\veps\|v_{xz}\|^2+\frac{\alf_2}{2}\|v_x(z=0)\|^2.
\end{split}
\end{equation}
 Finally, it follows from \eqref{e:vz.2}, \eqref{e:wvz.vzz} and
 \eqref{e:uvx.vzz.2} that, for $\veps>0$ chosen sufficiently small,
\begin{equation}
\label{e:vz.2c}
\begin{split}
 &\frac{d}{dt}\lf(\|v_z\|^2 +\alf_2\|v(z=0)\|^2 \rt)
 + \|\nbl v_z\|^2 + \alf_2\|v_x(z=0)\|^2\\
 \lsp & (\|u_x\|^2+\|v_x\|^2 + \|u\|^2\|u_z\|^2)\|v_z\|^2\\
 &+ \alf_2\|u_z\|\|u_{xz}\|\|v(z=0)\|^2 + \|u_z\|^4 + \|u\|^2\\
\lsp & (\|u_x\|^2+\|v_x\|^2 + \|u\|^2\|u_z\|^2+\|u_z\|\|u_{xz}\| )\\
 &\x(\|v_z\|^2+\alf_2\|v(z=0)\|^2) + \|u_z\|^4 + \|u\|^2.
\end{split}
\end{equation}
 Notice that, by \eqref{e:be.1}, for $v_{0,z}\in L^2$,
\[ \|v_0(z=0)\| \ls h^\hf\|v_{0,z}\|<\infty.\]
 By, Theorem~\ref{t:uz} (a), \eqref{e:vz.2c} and Gronwall lemma, we have local
 in time boundedness of $\|v_z\|$ on some interval $[0,t_0]$, with $t_0>0$.
 This yields global uniform boundedness and an absorbing set for $\|v_z\|$,
 since we have uniform boundedness and an absorbing set for $(u,v,\tht)$ in
 $V$ when considered in $[t_0,\infty)$ as a strong solution.

 Finally, with \eqref{e:vz.2c} we can justify as in our proof of part (a), that
 Theorem~\ref{t:uz} (b) is valid for a weak solution $(u,v,\tht)$ when
 $\Pz u_0$,$\Pz v_0\in L^2(D)$.

 Proof of Theorem~\ref{t:uz} (c) is similar to that for (b).

\end{prf}

 Notice that vertical regularity $(u_z,v_z,\tht_z)\in L^2$ of weak solutions
 of 2D and 3D viscous PE played a very prominent or even crucial role in almost
 all previous analytic works in dealing with solution regularity and uniqueness
 properties, for example, as shown in \eqref{e:sc}. However, to the contrary
 of this intuitive impression, we see next that horizontal regularity might
 actually force weak solutions to behave somewhat better, at least for the 2D
 problem. This is manifested in the following Theorem~\ref{t:ux}, which is our
 second main result of this section.

\begin{thrm}
\label{t:ux}
 Suppose $Q\in L^2(D)$, $(u_0, v_0,\tht_0)\in H$ and $(u,v,\tht)$ is a weak
 solution of \eqref{e:u}-\eqref{e:ic}. The following statements are valid:
\begin{enumerate}
 \item[{\em (a)}] If $\Px u_0\in L^2(D)$, then
 \[ u_x \in L^\infty(0, +\infty; L^2(D)) \cap L^2(0,\infty; H^1(D)). \]
 Moreover, there exists a bounded absorbing set for $u_x$ in $L^2(D)$.
 \item[{\em (b)}] If $(\Px u_0,\Px v_0)\in [L^2(D)]^2$, then
\[ (u_x,v_x) \in L^\infty(0, +\infty; [L^2(D)]^2)
  \cap L^2(0,\infty; [H^1(D)]^2).\]
 Moreover, there exists a bounded absorbing set for $(u_x,v_x)$ in $[L^2(D)]^2$.
 \item[{\em (c)}] If $(\Px u_0,\Px\tht_0)\in [L^2(D)]^2$, then
\[ (u_x,\tht_x) \in L^\infty(0, +\infty; [L^2(D)]^2)
  \cap L^2(0,\infty; [H^1(D)]^2).\]
 Moreover, there exists a bounded absorbing set for $(u_x,\tht_x)$ in
 $[L^2(D)]^2$.
 \item[{\em (d)}] If $(\Px u_0,\Pz v_0) \in [L^2(D)]^2$, then
\[ (u_x,v_z) \in L^\infty(0, +\infty; [L^2(D)]^2)
 \cap L^2(0,\infty; [H^1(D)]^2).\]
 Moreover, there exists a bounded absorbing set for $(u_x,v_z)$ in $[L^2(D)]^2$.
\item[{\em (e)}] If $(\Px u_0,\Pz\tht_0)\in [L^2(D)]^2$, then
\[ (u_x,\tht_z) \in L^\infty(0, +\infty; [L^2(D)]^2)
  \cap L^2(0,\infty; [H^1(D)]^2).\]
 Moreover, there exists a bounded absorbing set for $(u_x,\tht_z)$ in
 $[L^2(D)]^2$.
\end{enumerate}
 \end{thrm}

\begin{prf}

 Again, for simplicity of presentation, in the proof of Theorem~\ref{t:ux}, we
 will only provide the key estimates of $\|u_x\|$, $\|v_x\|$, $\|\tht_x\|$,
 $\|v_z\|$ and $\|\tht_z\|$ under the assumption that $(u,v,\tht)$ is a
 {\em strong} solution on $[0, \infty)$. The justification that these estimates
 are sufficient for a rigorous proof of Theorem~\ref{t:ux} is almost the same
 as the one we provided in our proof of Theorem~\ref{t:uz} (a). Thus, it is
 omitted for conciseness.

 We provide these key estimates in three steps.

 {\bf Step 1}. Estimate for $\|u_x\|_2$.\\
 Taking inner product of \eqref{e:u} with $-u_{xx}$ yields
\begin{equation}
\label{e:ux.2}
\begin{split}
 \hf\dt{\|u_x\|^2} +& \|\nbl u_x\|^2 + \alf_1\|u_x(z=0)\|^2\\
=& \lra{uu_x+wu_z+v +q_x+\int_z^0\tht_x(x,\zt,t)d\zt}{u_{xx}}.
\end{split}
\end{equation}
 The following computations, along with the boundary conditions
 \eqref{e:bc.u}-\eqref{e:bc.t} and the constraint \eqref{e:bc.u.2},
 are used in the derivation of \eqref{e:ux.2}:
\[ -\int_D u_t u_{xx} dxdz
 = -\int_{-h}^0\lf(u_tu_x\Big|_{x=0}^1 -\int_0^1 u_{xt}u_x\rt)
 = \hf\dt{\|u_x\|^2},\]
\begin{equation*}
\begin{split}
\int_D u_{xx}u_{zz} dxdz
 =& \int_{-h}^0 \lf( u_xu_{zz}\Big|_{x=0}^1 -\int_0^1 u_xu_{xzz} dx \rt) dz\\
 =& - \int_D u_xu_{xzz} dxdz\\
 =& -\int_0^1 \lf ( u_xu_{xz}\Big|_{z=-h}^0 - \int_{-h}^0|u_{xz}|^2 dz \rt)dx\\
 =& \alf_1\|u_x(z=0)\|^2 + \|u_{xz}\|^2,
\end{split}
\end{equation*}
\[ \int_D q_x u_{xx} dxdz=\int_0^1 q_x\lf(\int_{-h}^0u_{xx}dz \rt)dx
 =\int_0^1 q_x \Px\lf(\int_{-h}^0 u_x dz\rt) dx =0. \]
 Now, we estimate the two tri-linear terms on the right-hand side of
 \eqref{e:ux.2} as following:
\begin{equation*}
\begin{split}
\lra{uu_x}{u_{xx}}=& \int_{-h}^0 \int_0^1 u \Px\lf(\frac{u_x^2}{2}\rt) dxdz
= -\hf\int_D u_x^3\\
\ls& C\|u_x\|^2\|u_{xx}\|^\hf\|u_{xz}\|^\hf\\
\ls& C_\veps\|u_x\|^4+ \veps\|\nbl u_x\|^2.
\end{split}
\end{equation*}
\begin{equation*}
\begin{split}
\lra{wu_z}{u_{xx}} \ls& C\|w\|^\hf\|w_z\|^\hf\|u_z\|^\hf\|u_{zx}\|^\hf\|u_{xx}\|\\
\ls& C\|u_x\| \|u_z\|^\hf\|u_{zx}\|^\hf\|u_{xx}\|\\
\ls&  C_\veps\|u_z\|^2\|u_x\|^4+\veps\|\nbl u_x\|^2.
\end{split}
\end{equation*}
Therefore, it follows from \eqref{e:ux.2} that
\begin{equation*}
\begin{split}
\hf\dt{\|u_x\|^2} +& \|\nbl u_x\|^2  + \alf_1\|u_x(z=0)\|^2\\
\ls& 2\veps \|\nbl u_x\|^2 + C_\veps(1+\|u_z\|^2)\|u_x\|^4
 + (\|v\| + \|\tht_x\|)\|u_{xx}\|\\
\ls& 3\veps \|\nbl u_x\|^2
 + C_\veps(1+\|u_z\|^2)\|u_x\|^4 + C_\veps(\|v\|^2+\|\tht_x\|^2)
\end{split}
\end{equation*}
 Thus, if $\veps>0$ is chosen sufficiently small, then
\begin{equation}
\label{e:ux.2b}
\begin{split}
 \dt{\|u_x\|^2} +& \|\nbl u_x\|^2  + \alf_1\|u_x(z=0)\|^2\\
\lsp & (1+\|u_z\|^2)\|u_x\|^4 + \|v\|^2+\|\tht_x\|^2.
\end{split}
\end{equation}
 A local in time upper bound of $\|u_x\|$ can be obtained from \eqref{e:ux.2b}
 as follows. By \eqref{e:ux.2b},
\[\dt{\|u_x\|^2} \ls C (1+\|u_z\|^2 + \|v\|^2 + \|\tht_x\|^2)(\|u_x\|^2 +1)^2. \]
 Denote:
\[ y(t):=\|u_x\|^2+1, \quad g(t) := C (1+\|u_z\|^2 + \|v\|^2 + \|\tht_x\|^2).\]
 Then
\[ y'(t) \ls g(t) y^2(t). \]
 Notice that $y(t) \gs 1$, $g(t)\gs C(>0)$ and $g\in L^1(0,+\infty)$.
 Therefore,
\[ y(t) \ls \frac{y(0)}{ 1 - y(0)\int_0^t g(s) ds}, \]
 for $t\in(0, t_*)$ where $t_*$ is decided by the following equation:
\[  y(0) \int_0^{t_*} g(s) ds = 1. \]
 Thus, for $t\in(0, t_*)$
\[ \|u_x\|^2 +1 \ls
  \frac{ \|u_{0,x}\|^2 +1}{1 - (\|u_{0,x}\|^2+1)\int_0^t g(s)ds}. \]
 This finishes the proof of local in time boundedness of $\|u_x\|$.\\
 One the other hand, there exists $t_0\in(0, t_*)$, such that
\[ (u(t_0), v(t_0), \tht(t_0)) \in V. \]
 Therefore, by the result of uniform boundedness  of strong solutions (see
 \S 3.3 of \cite{p.t.z:09} and \cite{k.z:08}) and its uniqueness on
 $[t_0, +\infty)$, there exists a bounded absorbing set for $(u,v,\tht)$ in
 $V$ for $t\in[t_0,\infty)$, thus a bounded absorbing set for $u_x$ in $L^2$
 for $t\in[0,\infty)$. It also proves the uniform boundedness:
\[ u_x\in L^\infty(0, +\infty; L^2). \]
 Then, integrating \eqref{e:ux.2b} for $t$ from $0$ to $\infty$ proves
\[ \nbl u_x\in L^2(0,+\infty; L^2). \]
 This finishes proof of Theorem~\ref{t:ux} (a).
 
{\bf Step 2}. Estimate of $\|v_x\|_2$ and $\|\tht_x\|$.\\
 Similar to {\bf Step 1}, taking inner product of \eqref{e:v} with $-v_{xx}$
 yields:
\begin{equation}
\label{e:vx.2}
\begin{split}
 \hf\dt{\|v_x\|^2} + \|\nbl v_x\|^2 + \alf_2\|v_x(z=0)\|^2
= \lra{uv_x+wv_z-u}{v_{xx}}.
\end{split}
\end{equation}
 The tri-linear terms on the right-hand side of \eqref{e:vx.2} can be estimated as following:
\begin{equation*}
\begin{split}
\int_D uv_xv_{xx} =& \hf \int_D u\Px(v_x^2) dxdz\\
 =& -\hf\int_Du_x v_x^2 dxdz\\
\ls& C\|u_x\| \|v_x\|\|v_{xx}\|^\hf\|v_{xz}\|^\hf\\
\ls& C_\veps\|u_x\|^2\|v_x\|^2 + \veps\|\nbl v_x\|^2.
\end{split}
\end{equation*}
\begin{equation*}
\begin{split}
 \int_D w v_z v_{xx}
\ls& C\|w\|^\hf\|w_z\|^\hf\|v_z\|^\hf\|v_{xz}\|^\hf \|v_{xx}\|\\
\ls& C\|u_x\|\|v_z\|^\hf\|v_{xz}\|^\hf\|v_{xx}\|\\
\ls& C_\veps\|u_x\|^4\|v_z\|^2 + \veps \|\nbl v_x\|^2.
\end{split}
\end{equation*}
 Thus, if $\veps>0$ is chosen sufficiently small, then
\begin{equation}
\label{e:vx.2b}
\begin{split}
 \dt{\|v_x\|^2} +& \|\nbl v_x\|^2 + \alf_2\|v_x(z=0)\|^2\\
 \ls& C(\|u_x\|^2\|v_x\|^2 + \|u_x\|^4\|v_z\|^2 + \|u\|^2).
\end{split}
\end{equation}
 Similarly, we have
\begin{equation}
\label{e:tx.2b}
\begin{split}
 \dt{\|\tht_x\|^2} +& \|\nbl \tht_x\|^2 + \alf_0\|\tht_x(z=0)\|^2\\
 \ls& C(\|u_x\|^2\|\tht_x\|^2 + \|u_x\|^4\|\tht_z\|^2 + \|Q\|^2).
\end{split}
\end{equation}
 Notice the fact that $v_z, \tht_z\in L^2(0,\infty; L^2)$.
 Thus, as argued in {\bf Step 1},
 Theorem~\ref{t:ux} (b) follows from Theorem~\ref{t:ux} (a) and \eqref{e:vx.2b};
 Theorem~\ref{t:ux} (c) follows from Theorem~\ref{t:ux} (a) and \eqref{e:tx.2b}.

{\bf Step 3}. Estimate of $\|v_z\|_2$ and $\|\tht_z\|$.\\
 Recall \eqref{e:vz.2}. The two tri-linear terms on the right-hand side of \eqref{e:vz.2}  
 can be estimated in the following.\\
 First, we have 
\begin{equation}
\label{e:uvx.vzz}
\begin{split}
 \int_D uv_x v_{zz}
\ls& C\|u\|^\hf\|u_x\|^\hf \|v_x\|^\hf\|v_{xz}\|^\hf\|v_{zz}\|\\
\ls& C_\veps\|u\|^2\|u_x\|^2\|v_x\|^2 + \veps\|\nbl v_z\|^2
\end{split}
\end{equation}
 Next, the estimate \eqref{e:wvz.vzz} will be used again. Thus, if $\veps>0$ is
 chosen sufficiently small, then
\begin{equation}
\label{e:vz.2b}
\begin{split}
 \frac{d}{dt}\lf(\|v_z\|^2 +\alf_2\|v(z=0)\|^2 \rt)
 +& \|\nbl v_z\|^2 + \alf_2\|v_x(z=0)\|^2\\
 \lsp & \|u\|^2\|u_x\|^2\|v_x\|^2 + \|u_x\|^2\|v_z\|^2 + \|u\|^2.
\end{split}
\end{equation}
 Similarly, we have
\begin{equation}
\label{e:tz.2b}
\begin{split}
 \frac{d}{dt}\lf(\|\tht_z\|^2 +\alf_0\|\tht(z=0)\|^2 \rt)
 +& \|\nbl \tht_z\|^2 + \alf_2\|\tht_x(z=0)\|^2\\
 \lsp & \|u\|^2\|u_x\|^2\|\tht_x\|^2 + \|u_x\|^2\|\tht_z\|^2 + \|Q\|^2.
\end{split}
\end{equation}
 Now, Theorem~\ref{t:ux} (d) and (e) follow from \eqref{e:vz.2b} and
 \eqref{e:tz.2b} respectively as in {\bf Step 1} and {\bf Step 2}.

\end{prf}

\section{Uniqueness}
\label{s:un}
 In this section, we state and prove our last main result, Theorem~\ref{t:uniq},
 on uniqueness of weak solutions of 2D viscous PE \eqref{e:u}-\eqref{e:ic} when
 some initial partial regularity is assumed.
\begin{thrm}
\label{t:uniq}
 Suppose $Q\in L^2(D)$, $(u_0, v_0,\tht_0)\in H$ and $(u,v,\tht)$ is a weak
 solution of \eqref{e:u}-\eqref{e:ic}. Suppose further that one of the
 following initial regularity is valid:
\begin{align*}
  &(\Px u_0,\Px v_0,\Px\tht_0)\in (L^2(D))^3,\text{ or }\
   (\Px u_0,\Px v_0,\Pz\tht_0)\in (L^2(D))^3,\\
\text{ or }\ &(\Px u_0,\Pz v_0,\Px\tht_0)\in (L^2(D))^3,\text{ or }\
   (\Px u_0,\Pz v_0,\Pz\tht_0)\in (L^2(D))^3,\\
\text{ or }\  &(\Pz u_0,\Pz v_0,\Pz\tht_0)\in (L^2(D))^3.
\end{align*}
 Then, the following are valid:
\begin{enumerate}
 \item[{\em (a)}] The norm for the corresponding solution regularity is
 uniformed bounded for all time $t\gs0$ and an absorbing set exists for the
 norm of the corresponding solution regularity.
 \item[{\em (b)}] The weak solution is unique.
\end{enumerate}
\end{thrm}

 \begin{prf}

 The claim (a) for solution regularity is an immediate consequence of
 Theorem~\ref{t:uz} and Theorem~\ref{t:ux}. Moreover, assuming any one of the
 above five initial conditions, we have
\begin{equation}
\label{e:xz}
 (u_{xz}, v_{xz}, \tht_{xz}) \in \lf[ L^2(0,\infty; L^2(D)) \rt]^3.
\end{equation}
 Notice that
\begin{equation*}
\begin{split}
 \|u_z\|_{L^\infty_x(L^2_z)}
=& \lf\|\lf(\int_{-h}^0|u_z|^2\rt)^\hf\rt\|_{L^\infty_x}\\
\ls & \text{ (Minkowski inequality) } \\
\ls & \lf(\int_{-h}^0 \|u_z\|_{L^\infty_x}^2\rt)^\hf\\
\ls & \text{ (Agmon's inequality) } \\
\ls & C\lf(\int_{-h}^0 \|u_z\|_{L^2_x}\|u_{xz}\|_{L^2_x}\rt)^\hf\\
\ls & C\|u_z\|^\hf\|u_{xz}\|^\hf.
\end{split}
\end{equation*}
 Therefore, for any $T>0$,
\begin{equation}
\label{e:uz.2}
\begin{split}
 \int_0^T \|u_z\|_{L^\infty_x(L^2_z)}^2 dt
 \ls& C\int_0^T \|u_z\|\|u_{xz}\| dt\\
 \ls& C\lf(\int_0^T \|u_z\|^2\ dt\rt)^\hf
 \lf(\int_0^T \|u_{xz}\|^2\ dt\rt)^\hf.
\end{split}
\end{equation}
 Notice that for any weak solution $(u,v,\tht)$,
\[ (u_z,v_z,\tht_z)\in L^2(0,T; L^2(D)). \]
 Thus, by \eqref{e:xz} and \eqref{e:uz.2}, we have
\[ u_z \in L^2(0,T; L^\infty_x(L^2_z)). \]
 Similarly,
\[ (u_z, v_z, \tht_z) \in \lf[L^2(0,T; L^\infty_x(L^2_z))\rt]^3. \]
 Thus, Theorem~\ref{t:uniq} (b) is proved by Theorem~\ref{t:sc}.

 \end{prf}


\end{document}